\documentclass{amsart}
  \usepackage{amscd,amssymb,epsfig}
   \usepackage{epic,eepic}

\input{amssym.def}
\input{amssym}
\usepackage{amscd}

                     \numberwithin{equation}{subsection}

                     \newtheorem{propo}{Proposition}[subsection]
                     \newtheorem{corol}[propo]{Corollary}
                     \newtheorem{theor}[propo]{Theorem}
                     \newtheorem{lemma}[propo]{Lemma}
                     \theoremstyle{definition}

                     \theoremstyle{remark}

                     \newcommand{\ZZ}{\mathbb{Z}}
                     \newcommand{\RR}{\mathbb{R}}

                     \newcommand{\Hom}{\operatorname{Hom}}
		     \newcommand{\Ker}{\operatorname{Ker}}

\newcommand{\rk}{\operatorname{rank}}
\newcommand{\Ann}{\operatorname{Ann}}
\newcommand{\card}{\operatorname{card}}
\newcommand{\eee}{{S}}
\newcommand{\mmm}{{C}}

		    \newcommand{\Ima}{\operatorname{Im}}
		      \newcommand{\inc}{\operatorname{in}}
		      
                     \newcommand{\id}{\operatorname{id}}

		     \newcommand{\sign} {\operatorname {sign}}
		    \newcommand{\modu}{\operatorname{mod}}

\begin{document}
      \title{Cobordism of knots on surfaces}
                     \author[Vladimir Turaev]{Vladimir Turaev}
                     \address{%
              IRMA, Universit\'e Louis  Pasteur - C.N.R.S., \newline
\indent  7 rue Ren\'e Descartes \newline
                     \indent F-67084 Strasbourg \newline
                     \indent France \newline
                     \indent  and \newline
                       \indent  Department of Mathematics \newline
                         \indent  Indiana University \newline
                           \indent   Rawles Hall,  831 East 3rd St  \newline
    \indent  Bloomington, IN 47405 \newline
                             \indent  USA \newline  }
                     \begin{abstract} 	We introduce a relation of cobordism for  knots in thickened surfaces and  study cobordism invariants of such knots.  

 

                     \end{abstract}

                     \maketitle

%%%\centerline {\bf Contents}
%%% \vskip1truecm

%%%\noindent  {\bf  1. Introduction}

                  \section*{Introduction} 
		  
	This paper is concerned with knots on oriented surfaces, that is with   knots in 3-manifolds of type (an oriented  surface) $\times \,\RR$.  The class of such knots  appears to be a natural intermediate between the class of classical knots   in $\RR^3=\RR^2\times \RR$  and  the class of   knots in arbitrary   3-manifolds.  The    diagrammatical methods used for    knots in $\RR^3$     extend  to knots 
	on surfaces.  On the other hand, one encounters  new phenomena absent in the classical case. The study of knots on surfaces has a long history, see for instance \cite{tu1}, \cite{crr}, \cite{fi}, \cite{cks}. This area has been especially active in the context of the Kauffman theory of virtual knots, see \cite{Ka}, \cite{km}, \cite{ku}. For other approaches, see  \cite{aps}, \cite{ct}, \cite{ft}.
	
The principal aim of this paper is to introduce and to study an equivalence relation of cobordism  for knots on surfaces.  Briefly speaking,  a knot   $K_1$ on a surface~$\Sigma_1$ is cobordant to a knot   $K_2$ on a surface $\Sigma_2$  if there is an oriented 3-manifold $M$ with $ \Sigma_1 \amalg (-\Sigma_2)\subset \partial M$ such that  the 1-manifold  $K_1\amalg (-K_2)\subset \partial M \times \RR$ bounds an  annulus in $M \times \RR$. We    introduce several  non-trivial cobordism invariants of knots. This includes geometric invariants (slice genera) and algebraic invariants:  polynomials $u_+, u_-$, algebraic genera, and a so-called graded matrix.    The   algebraic invariants  of a knot   on a surface $\Sigma$  are derived from an arbitrary diagram $D$ of this knot on $\Sigma$. The key idea is that every crossing   of $D$ splits the underlying loop of $D$ into two \lq \lq half-loops", and we can select one of them in a canonical way.  This gives a finite family of loops on $\Sigma$. The homological    intersection numbers   of these loops with each other and with   the underlying loop of $D$ form      a skew-symmetric square matrix over $\ZZ$. It  will be the main source of   our   invariants.  
	
 	It is natural to compare the 	  knot cobordism defined above with the standard   concordance for    knots in $\RR^3$, see \cite{ch} for a recent survey. It is obvious that concordant classical knots are cobordant
	in our sense. Besides this fact, essentially nothing seems to be  known about our relation of  cobordism  for classical knots. To the best of the author's knowledge, it is possible that all classical knots are cobordant to each other. It is also conceivable that two classical knots are cobordant if and only if they are concordant in the  standard sense. All knot invariants defined in this paper are trivial for the classical knots.

	One way to obtain invariants of  a knot $K\subset \Sigma \times \RR $ is to project $K$ to   $\Sigma$ and to study homotopy invariants of the resulting     loop  $\underline K$ of $K$.  Loops on surfaces and their cobordisms were     studied   in \cite{tu2}.   The  present paper  refines the methods of \cite{tu2} and lifts them to knots.  Preliminary knowledge of \cite{tu2} can be useful  to the reader  but  is not required.

 The content of the paper is as follows.  As a warm up, we introduce in   Section~\ref{1} the knot polynomials 
 $u_+$ and $u_-$. In Section \ref{2} we discuss cobordisms  of knots and show that $u_+$ and $ u_-$ are cobordism invariants. In Section \ref{3}  we  study graded matrices  and   introduce  graded matrices of knots.   In Section \ref{5}      we   analyze  an equivalence relation of cobordism on the class of graded matrices.  The relationship with  knot cobordisms  is discussed  in Section \ref{6}.
 Miscellaneous remarks and open questions   are collected in Section \ref{7}.
 
 Throughout the paper we work in the smooth category though all the results can be reformulated in PL and topological categories.  All surfaces and 3-manifolds in this paper are oriented.

		       \section{Polynomials $u_+$ and $u_-$}\label{1}
                    
		         \subsection{Knots and knot  diagrams}\label{p1}   
		         By a  {\it  knot}  	we  mean a pair consisting of an (oriented) surface $\Sigma$ and an oriented embedded circle $K\subset \Sigma \times \RR $ disjoint from $\partial \Sigma \times \RR$.  We   say that $K$ is a {\it knot   on} $\Sigma$. 
This knot will be   denoted $K\subset \Sigma \times \RR$ or $(\Sigma, K)$, or simply $K$. Note that $\Sigma$ may be non-compact. 

		         Two knots $ K_1, K_2\subset \Sigma \times \RR $   on a surface $\Sigma$ are {\it isotopic} if~$K_1$ can be deformed into $K_2$ in the class of knots in $   \Sigma \times \RR $ (keeping the orientation of the knot).
 Two knots $K_1 \subset \Sigma_1 \times \RR $ and $  K_2\subset \Sigma_2 \times \RR $ are {\it diffeomorphic} if 
 		         there is an orientation preserving diffeomorphism $f:\Sigma_1 \to \Sigma_2$ 
          such that the knot $(f\times \id_\RR)(K_1)\subset \Sigma_2 \times \RR $ is isotopic to $K_2$.
         Isotopic knots are necessarily  diffeomorphic.
		         
		         		         A   {\it {knot diagram}} on   a  surface  $\Sigma$ is an oriented   closed curve on $\Sigma - \partial \Sigma$ with only double transversal crossings such that at
each  double point,  one of the branches of the curve passing through this point is distinguished. The distinguished branch is
said to be {\it overgoing} while the second branch passing through the same double point is  {\it undergoing}.  A knot diagram $D$ on
$\Sigma =\Sigma\times
\{0\}$ determines a  knot in  
$\Sigma\times
\RR$ by    pushing the overgoing branches of $D$  into  $\Sigma\times (0,\infty) $.  It is clear that every knot in 
$   \Sigma \times \RR $ is isotopic to a knot presented by a knot diagram on $\Sigma$. Two knot diagrams on $\Sigma$ present isotopic knots if and only if these diagrams are related by a finite sequence of standard Reidemeister moves (briefly, R-moves) applied in small disks in $\Sigma$  while  keeping the rest of the diagram.

 \subsection{Halfs of   knot diagrams}\label{p2} 	Let $D$ be a knot diagram on a  surface $\Sigma$. 
  The set of double points of $D$ will be denoted 	$ \Join\!\!  (D)$. 
 With each double point $x\in \, \,\Join\! \!  (D)$, we associate a loop $D_x$ on  
   $ \Sigma   $      as follows. Let    $A$ and $B$ be the     branches of $D$ passing through $x$, where the notation is chosen  so that the pair (a   positive  tangent vector of $A$ at $x$, a   positive  tangent vector of $B$ at $x$) is a positive basis in the tangent space of $\Sigma$ at $x$. Denote by $D_x$   the loop on $ \Sigma$    starting at~$x$, going  along $A$ in the positive direction, and then going  along the underlying loop of~$D$ 
  until the first return to $x$. We    call $D_x$ the {\it {distinguished half}} of $D$ at $x$.  Note that the definition of $D_x$ does not use the over/under data at $x$. To keep record of this data, we define  the {\it sign} $\sign(x)$ of   $x $ to be   $ +1$ if $A$ is   overgoing  and $-1$ if $B$ is   overgoing. The set $\Join\!\!  (D)$ splits as a  union of two disjoint subsets $$\Join_+\!\!  (D)=\{x\in \,\,\Join\!\!  (D)\,\vert \, \sign (x)=+1\}\,\,\, \text{and} \,\,\, \Join_-\!\! (D)=\{x\in \,\,\Join\!\!  (D)\,\vert \, \sign (x)=-1\}.$$ 
		         
 \subsection{Definition  of  $u_\pm$}\label{1.3} Let $K\subset 	\Sigma \times \RR$  be a knot presented by a diagram $D$ on~$\Sigma$. 	Let $[D]\in H_1(\Sigma)=H_1(\Sigma; \ZZ)$ be the homology class of the (oriented) closed curve underlying $D$.
 For a  double point 	$x \in \, \,\Join\!\!  (D)$,  we can   consider the homology class $[D_x]\in 	H_1(\Sigma)$
 of the loop $D_x$    and the homological intersection number $     [D_x] \cdot [D] \in \ZZ$.
 For shorteness, we shall write $ D_x  \cdot  D $ for $ [D_x] \cdot [D]$.
 
We define two  one-variable polynomials $u_+(K), u_-(K) \in \ZZ[t]$ by
 $$u_\pm(K)=\sum_{x\in \,\Join \!   (D), \, D_x  \cdot  D\neq 0,\, \sign (D_x  \cdot  D)=\pm \sign (x) } \sign (x) \,t^{ \vert \,D_x  \cdot  D \,\vert },$$
 where $ \sign (D_x  \cdot  D)\in \{+1, -1\}$ and $\vert \,D_x  \cdot  D \,\vert$ are respectively the sign and the absolute value of the non-zero integer $  \,D_x  \cdot  D $. In other words, 
 $$u_+(K)=\sum_{x\in \,\Join_+\!  (D), \, D_x  \cdot  D>0} t^{D_x  \cdot  D}
 - \sum_{x\in \,\Join_-\! (D), \,D_x  \cdot  D<0} t^{- (D_x  \cdot  D)} ,$$
$$u_-(K)=   \sum_{x\in \,\Join_+\! (D), \,D_x  \cdot  D<0} t^{- (D_x  \cdot  D)} -\sum_{x\in \,\Join_-\!   (D), \, D_x  \cdot  D>0} t^{D_x  \cdot  D}
 .$$	
 Both   $u_+(K)$ and $u_-(K)$ are isotopy invariants of $K$. To see this, it is enough to check that   $u_{\pm}(K)$ is preserved under the R-moves on $D$. The invariance under the first 	R-move creating a new crossing $x_0$ follows from the fact that $[D_{x_0}]=0$ or      $[D_{x_0}]=[D]$; in both cases  $ D_{x_0}  \cdot  D =0$.      
 The invariance under the second 	R-move creating two new crossings $x_1, x_2$ follows from the fact that $\sign (x_1)=-\sign (x_2)$ and $[D_{x_1}] + [D_{x_2}]=[D]$ so that $ D_{x_1}  \cdot  D =-(D_{x_2} \cdot  D)$.  If $ D_{x_1}  \cdot  D =  D_{x_2}  \cdot D =0$, then   $x_1,x_2$ do not contribute to $u_{\pm}(K)$.    If $ D_{x_1}  \cdot  D =-(D_{x_2}  \cdot  D)\neq 0$, then the   contributions of $x_1,x_2$ cancel each other. The invariance under the third 	R-move follows from the fact that the set of pairs $\{([D_x], \sign (x))\}_x$ is the same before and after the move. 
Thus, the polynomials $u_+, u_- $ are isotopy invariants of knots.  It is also clear  that they are diffeomorphism invariants of knots.
 
 \subsection{Properties  of  $u_\pm$}\label{1.4} We  point out several properties of  the polynomials $u_{\pm}(K)$ of a knot $K$ on a surface $\Sigma$.
  It is clear that   $u_+(K)$ and $u_-(K)$ have zero free terms.
 The difference $u(K)=u_+(K) -u_-(K)$ is computed from a diagram $D$ of $K$ by
 $$u(K)=\sum_{x\in \,\Join\!   (D), \, D_x  \cdot  D \neq 0} \sign (D_x  \cdot  D)\, t^{\vert D_x  \cdot  D\vert }
  \,.$$
 It follows from this formula that $u(K)$ does not depend on the   over/under data in   $D$.  Therefore, the 
  polynomial  $u=u(K)\in \ZZ[t] $ is a homotopy invariant of $K$. This 
 polynomial   is the   invariant of the   underlying loop of~$D$ on~$\Sigma$    introduced in \cite[Section 3]{tu2}. One property of $u$ established there computes the value of its derivative $u'=   {du}/{dt}$ at $t=1$, namely $u'(1)=0$. We conclude that  
\begin{equation}\label{rede} u'_+(K)(1)=  u'_-(K)(1).\end{equation}
 The reader can prove this equality directly as an exercise. 
 
 We  denote by  $\overline K$   the same knot $K\subset \Sigma \times \RR$ with opposite orientation on~$\Sigma$ (keeping the one on~$K$). Similarly,  
 $-K$ denotes   the same knot $K\subset \Sigma \times \RR$ with opposite orientation on $K$ (keeping the one on~$\Sigma$).  It is easy to check that 
 $u_\pm (\overline K)=- u_\mp(K)$ and 
 $u_\pm (-K)= u_\mp (K)$. If two knots   $K_1, K_2 \subset \Sigma \times \RR$ are presented by disjoint  diagrams $D_1, D_2$ on $\Sigma$   and a knot  $K_1\#K_2\subset \Sigma \times \RR$ is presented by the connected sum of $D_1, D_2$ along an   embedded arc on $\Sigma$ relating a point on $D_1$ to a point on $D_2$ and disjoint from $D_1, D_2$ otherwise, then 
 $u_{\pm} (K_1\#K_2)=u_{\pm} (K_1) + u_{\pm} (K_2)$. 
 
  \subsection{Examples}\label{exal} 1. For any knot   $K\subset \RR^3=\RR^2\times \RR$, we have $u_+(K)=u_-(K)=0$. Indeed, the homological intersection number of any two loops on $\RR^2$ is zero.
  
  2.  Let  $D$ be the knot diagram on the torus $S^1\times S^1$   obtained from the standard diagram of a trefoil on $S^2$ by replacing  the over-going branch at one crossing   with a branch going along a  1-handle attached to $S^2$. The diagram $D$ has only two double points $x, y$. A direct computation shows that for   appropriate orientations on $D$ and on the torus, $\sign(x)=\sign (y)= D_x  \cdot  D =- D_y   \cdot  D = 1$. Then for   the knot   $K\subset S^1\times S^1 \times \RR$ presented by $D$, we have $u_+(K)= u_-(K)= t$.  Note that $u_+(K)= u_-(K)$ for an arbitrary knot 
  $K\subset S^1\times S^1 \times \RR$;  this follows from the fact that any loop on the torus may be deformed inside an annulus and the intersection number of any two loops on an annulus is  zero.
  
 \subsection{Realization}  We say that a pair of  polynomials  $p_+, p_-\in \ZZ[t]$ is {\it {realized}} by a knot
 $K$ if $u_+(K)=p_+$ and $u_-(K)=p_-$. If a pair $(p_+, p_-)$ is realized by a knot, then it can be realized by a knot on a closed oriented surface. This follows from the facts that the computation of  $u_+, u_-$ proceeds entirely inside   a neighborhood of a knot diagram in the ambient surface and that any  knot diagram on a surface has a neighborhood which embeds into a closed surface.
  
 The following theorem solves the realization problem for  $u_+ , u_- $.

 \begin{theor}\label{th1}
                     A pair of  polynomials  $p_+, p_-\in \ZZ[t]$   can be realized by a knot    if and only if $p_+(0)= p_-(0)=0$   and $p'_+(1)=p'_-(1) $. 
                     \end{theor}
                     \begin{proof} We only need  to prove the sufficiency. Set $u=p_+ -p_-$. The conditions on $p_+, p_-$ imply that  
                     $u(0)= 0$  and  $u'(1)=0$. It follows from  
                     \cite[Theorem~3.4.1]{tu2}, that there is a knot diagram  $D$ on    an oriented surface $\Sigma$ presenting a knot   $K_0\subset \Sigma \times \RR$ with $u(K_0)=u_+(K_0)-u_-(K_0)
                      =u$.   The knot   $K_0 $ and the knot  $K_0 \# (S^1\times x)$ on  $\Sigma \# (S^1 \times S^1)$,  where $x\in S^1$,  have the same   $u_\pm$. Therefore, replacing if necessary $K_0 $ by   $K_0 \# (S^1\times x)$,  we can assume that  there is a loop $\alpha$ on $\Sigma$    such that $ \alpha \cdot  D =1$. We claim that      for  any 
                     knot $K'\subset \Sigma \times \RR$ homotopic to $K_0$ and any non-zero integer $N$, there is a diagram $D'$ of $K'$ on          $\Sigma$ and a double point $x\in \,\Join \!\! (D')$ such that  $$\sign (x)=\sign (N) \quad \text {and} \quad D'_x\cdot D'=N.$$  This claim   implies the theorem. Indeed, exchanging the over/under-crossings at such   $x\in \,\Join \!\! (D')$, we obtain a new knot diagram on $\Sigma$ presenting a knot $K''\subset \Sigma \times \RR$ homotopic to $K_0$. Clearly, $u_+(K'')=u_+(K') - \sign (N)\,  t^{\vert N\vert}$. Starting from $K_0$ and applying this procedure recursively, we can obtain a knot $K\subset \Sigma \times \RR$ homotopic to~$K_0$ and such that $u_+(K)=p_+$.  Then  $$u_-(K)=u_+(K)-u(K)= p_+- u(K_0)= p_+-u=p_-.$$
  To prove the claim above, we take a small arc on $D$ and push it with a finger    along a generic loop
  on $\Sigma$   homologous to $\alpha^N$ and passing everywhere over (or everywhere under) $D$.   When the finger comes back to $D$ near the original arc we  obtain two new self-crossings; one of them will be the required $x$. We leave the details to the reader.                 
\end{proof}

\subsection{Coverings of knots and higher polynomials} The polynomials $u_\pm$ give  rise to a   family of polynomial invariants of knots
numerated by finite sequences of  positive integers. We first define the $m$-th covering  of a knot $(\Sigma, K)$ for any integer $m\geq 1$.   Projecting the homology class $[K]\in 
H_1(\Sigma\times \RR; \ZZ )=H_1(\Sigma;   \ZZ)$ to $H_1(\Sigma; \ZZ/m\ZZ)$ and applying the  Poincar\'e duality $H_1(\Sigma;   \ZZ/m\ZZ) \cong H^1(\Sigma ; \ZZ/m\ZZ)$, we obtain a cohomology class in $ H^1(\Sigma ; \ZZ/m\ZZ)$. This class determines an  $m$-fold covering $  \Sigma^{(m)} \to \Sigma$. It is clear that    $K$ lifts to 
a knot $  K^{(m)} \subset  \Sigma^{(m)} \times \RR$. This knot    is determined by   $m$ and the diffeomorphism type of $(\Sigma, K)$   uniquely up to diffeomorphism. 
Iterating this construction, we  define for any finite sequence of   integers $m_1,\ldots, m_k \geq 1$, a knot
$$ K^{(m_1,\ldots,m_k)}= (\ldots(   K^{(m_1)})^{(m_2)}\ldots)^{(m_k)}.$$ 
Set 
$$u^{m_1,\ldots, m_k }_\pm ( K)= u_\pm (  K^{(m_1,\ldots,m_k)})\in \ZZ[t].$$
By the results above,  this polynomial  is a diffeomorphism invariant of $(\Sigma, K)$.
Note that $u^1_+(K)=u_+ (K)$ and $u^1_-(K)=u_- (K)$.

			\section{Cobordism  and sliceness of knots}\label{2}

\subsection{Cobordism of knots}\label{sn:nbee}  
 Two   knots  $(\Sigma_1, K_1)$ and $(\Sigma_2, K_2)$ are   {\it cobordant}  if there is an oriented 3-manifold $M$ and an embedding  $ \Sigma_1 \amalg (-\Sigma_2) \hookrightarrow \partial M$ such that  the image of $K_1\amalg (-K_2)$ under the induced embedding
 $( \Sigma_1\times \RR) \amalg ( \Sigma_2\times \RR)  \hookrightarrow  \partial M \times \RR$ bounds an embedded oriented annulus in $M \times \RR$. It is understood that the orientation of $M$ (resp.\ of the annulus) induces the given orientation on $\Sigma_1$ (resp.\ on $K_1$) and the orientation opposite to the given one on $\Sigma_2$ (resp.\ on $K_2$).    A standard argument shows that   cobordism of knots is an equivalence relation. 
%% If   $ K_1$ is cobordant to  $K_2$, then $\overline K_1$ is cobordant to $\overline K_2$ and %%$K_1^-$ is cobordant to $K_2^-$.  

It follows from the definitions that diffeomorphic knots are cobordant. If $\Sigma'$ is a subsurface of a surface $\Sigma$, then any knot   in $   \Sigma'\times \RR$ is cobordant to the same knot viewed 
in $\Sigma \times \RR$. (To see this,  take $M=\Sigma \times [0,1]$ and consider the embeddings $\Sigma=\Sigma \times 0 \hookrightarrow \partial M$  and $\Sigma'=\Sigma' \times 1 \hookrightarrow \partial M$.)  In the language of \cite{cks}, these facts may be reformulated by saying that stably equivalent knots are cobordant.

 \begin{lemma}\label{le} If two knots    $  K_1 $  and $  K_2 $ are cobordant, then the knots $  K_1^{(m)}$ and $   K_2^{(m)}$ are  cobordant  
 for any $m\geq 1$.  \end{lemma}

\begin{proof} Let $M$ be an oriented 3-manifold as in the definition of cobordism, so  that    the  1-manifold    $K_1\amalg (-K_2)  \subset \partial M \times \RR$ bounds an embedded oriented annulus   $F\subset M \times \RR$.  Let $p:M\times \RR \to M$ be the projection.  
Replacing, if necessary,  $M$ by a closed neighborhood of   $p(F)$, we can 
assume that $M$ is compact.  By the Poincar\'e duality,  there is a cohomology class $a\in H^1(M;\ZZ)$ dual to the homology class of $p(F)$ in $  H_2(M,\partial M; \ZZ)$.  Let $  M^{(m)}\to M$ be the $m$-fold covering determined by $a\, (\modu m)\in H^1(M; \ZZ/m\ZZ)$. 
The annulus $F$ lifts to an embedded annulus in $  M^{(m)}\times \RR$; the latter yields a cobordism between $K_1^{(m)}$ and $K_2^{(m)} $.
\end{proof}

 \subsection{Slice knots}\label{fbm}    A   knot $K\subset \Sigma\times \RR$  cobordant to a trivial knot in $\RR^3=\RR^2\times \RR$ is said to be   {\it slice}.  It is easy to check that a  knot $K\subset \Sigma\times \RR$ is slice  if there is an oriented 3-manifold $M$ and an embedding  $ \Sigma \subset \partial M$ such that  the image of $K $ under the induced embedding
 $  \Sigma \times \RR  \subset \partial M \times \RR$ bounds an embedded disk in $M \times \RR$.

One should be careful in applying this terminology to   knots in $\RR^3$. A  knot  in $\RR^3$ that is slice in the usual sense of the word is   slice in our sense, but the converse is unknown and possibly untrue. In the sequel we use the term \lq\lq slice" exclusively in the sense defined above.

It is clear that    knots  cobordant  to a slice knot are 
slice.     By Lemma \ref{le}, the coverings of   slice knots are slice. 
  Note also that the usual notion of ribbon knots in $\RR^3$ directly extends to knots on surfaces  and yields a wide class of slice knots.

 \subsection{Polynomials of cobordant   knots}\label{fkkkm} We state the main result of this section.

 \begin{theor}\label{th:0388} The polynomials $u_+$ and $u_-$ are cobordism invariants of knots.   \end{theor}

This theorem follows from the identity  $u_\pm (-\overline K)=-u_\pm(K)$ and the following  more general lemma.

	\begin{lemma}\label{l:nblddddd}     Let $M$ be an oriented  3-manifold and $K_1,\ldots, K_r $ be  knots in $\partial M \times \RR$ whose projections   to  $\partial M$ are pairwise disjoint.  If the oriented 1-manifold  $\cup_{i=1}^r K_i \subset  \partial M \times \RR$ bounds in $   M\times \RR$    an embedded    compact   oriented  surface of genus $0$, then
\begin{equation}\label{osn} \sum_{i=1}^r u_\pm (K_i)=0. \end{equation}
\end{lemma} 
                     \begin{proof}    We need some terminology  concerning  maps   from     surfaces to $M$.        We say that  a map   $\omega$ from a surface $F$ to $M$    is    {\it
generic} if   

(i)   $\card  (\omega^{-1}(x))\leq 3$ for all $x\in M$;

 (ii)  $\omega^{-1} (\partial M)=\partial F$ and the restriction of $\omega$ to   $\partial F$ is an immersion into  $\partial M$ with only double transversal  crossings; any point of $ \omega(\partial F )$ has  a
neighborhood
$V\subset M$ such that the pair $(V, V\cap \omega(F))$ is diffeomorphic to    $(\RR^2, \RR  \times 0)  \times \RR_+
$  or to 
$(\RR^2,  \RR  \times 0  
  \cup 0\times
\RR )\times \RR_+$, where  $  \RR_+=\{r\in \RR\, \vert \, r\geq 0\}$;

(iii) there is a finite set $B\subset   F^\circ=F-\partial F$ such that the restriction of $\omega$ to  
$F^\circ-B $ is   an immersion into   $M^\circ=M-\partial M$ and any   distinct points $a,a'\in F^\circ-B $ with $\omega(a)=\omega (a')$ have disjoint neighborhoods in $F$ whose images under $\omega$ are transversal; 
                   
   (iv) for any $b\in B$,      there is a closed 3-ball $D\subset  M^\circ$
such that   $\omega(F)\cap \partial D$ is  a figure 8 curve   (a closed curve with one transversal self-crossing) and 
$\omega(F)\cap D$ is the cone over this curve with cone point $\omega(b)\in  D-\partial D$.

The points of the set $B=B(\omega)$ are called the {\it   branch points} of   $\omega$.  Conditions (i) and (iii)   imply that          any point of $ \omega(F^\circ-B )$ has  a neighborhood $V\subset M^\circ$ such that the pair $(V, V\cap \omega(F))$ is
diffeomorphic to the pair  ($\RR^3$, a union of $i$ coordinate planes) for $i=1,2, 3$. For more on this, see \cite[Chapter 7]{Kam}.

 Let now $F\subset M\times \RR$ be an embedded    compact      surface of genus $0$ 
 bounded by $\cup_{i=1}^r K_i$. We can slightly deform $F$ in the class of   embedded surfaces with boundary on $\partial M \times \RR$ to ensure  that  $ F\cap (\partial M \times \RR)=\partial F$  and  the restriction to $F$ of the projection $M\times \RR \to M$ is a generic map. We can still assume  that the components of $\partial F$ in $\partial M  \times \RR$ project to disjoint
 subsets of $\partial M$. These components are isotopic to the knots $K_1,\ldots, K_r$ and have the same polynomials $u_\pm$. Thus,  to prove the lemma, we can use the same symbols $K_1,\ldots, K_r$ for the components of~$\partial F$.  The   generic map  from $F$ to $M$ obtained as the restriction   of the projection $M\times \RR \to M$ will be  denoted $\omega$. Set  $B=B(\omega)\subset F^\circ$ and  $$ \mmm=B \cup  \{a\in F\,\vert\,   \card\, ( \omega^{-1}
(\omega(a)))
\geq 2\}\subset F.$$   It is clear from the conditions (i) -- (iv) above that  
$ \mmm$ consists of a finite number of   immersed  circles     and       intervals  in $F$ meeting $\partial F$ exactly at the endpoints of the intervals.  These endpoints   are the  preimages of the double points of
$\omega\vert_{\partial F}:\partial F \to  \partial M$.  The set of these endpoints coincides with $ \mmm \cap \partial F$ and is denoted $\partial  \mmm$. The immersed circles and intervals forming  $ \mmm$  have only double transversal crossings and self-crossings all lying in $F^\circ$.   The set of these crossings and self-crossings is denoted 
$\eee$; it    consists precisely of the preimages of the triple points of
$\omega$.      The   subsets $B , \eee,\partial \mmm $ of $\mmm$ are finite and pairwise disjoint.

Let  $\widetilde \mmm$ be  an abstract  1-dimensional manifold parametrizing $ \mmm$. The   projection $p: \widetilde \mmm
\to  \mmm$ is 2-to-1 over  $\eee$ and   1-to-1 over $ \mmm-  \eee$.   
We shall identify $\widetilde \mmm-p^{-1} (\eee) $ with  $ \mmm - \eee$ via $p$. 

For any  point $a\in \mmm -(\eee\cup B)$
there is exactly one   other point
$a'  \in \mmm-(\eee\cup B)$ such that $\omega(a)=\omega(a' )$.  The correspondence $a    \mapsto a' $  on 
$$\mmm -(\eee\cup B)=\widetilde \mmm - p^{-1} (\eee\cup B)$$ extends by continuity
to an involution 
   on
$\widetilde \mmm$, denoted $\tau$. The set   of   fixed
points  of $\tau$  is    $B$.  It is clear that $\omega\, p\,\tau=\omega\, p: \widetilde \mmm\to M$.  For   $a\in \partial
\mmm=\partial 
\widetilde \mmm$, the  point
$ \tau(a)  $ is the unique other point  of  $ \partial \mmm $   such that $\omega(a)=\omega(\tau (a) )$.  
Applying $\omega$, we can identify the quotient set $\partial \mmm/\tau$   with the set $\Join
\! \!(\omega\vert_{\partial F})\subset \partial M$ of  the self-crossings
of $\omega\vert_{\partial F}$.

We   define an  involution  $\mu:\partial \mmm\to   \partial \mmm$.   For   $a\in \partial \mmm=\partial 
\widetilde \mmm$,  let $I_a \subset  \widetilde \mmm$ denote the  
  interval    adjacent to $a$ and  let  $\mu(a)\in
\partial \mmm$     be its   endpoint   distinct from~$a$.  Clearly, $\mu^2=\id$ and $\mu$ is  fixed-point-free.   Note that $\mu$ commutes with $\tau\vert_{\partial \mmm}$.
Indeed, if  $
\tau(I_a)=I_a$, then   $\tau$ exchanges the endpoints of
$I_a$ so that   $ \tau=\mu $ on $\partial I_a\subset \partial \mmm$.  (In this case $\tau$ must have a unique fixed point on $I_a$,  so that    
$I_a$ contains   a unique branch point of $\omega$.)  If $\tau(I_a)\neq I_a$,  then  $\tau(I_a)$ has the endpoints $\tau(a), 
\tau(\mu(a) )$ so that $\mu(\tau(a))=\tau (\mu(a))$.    Since  $\mu$ and $\tau \vert_{\partial \mmm}$ commute,   $\mu$
induces an involution on   $\partial \mmm/\tau=\,  \Join
\! \!(\omega\vert_{\partial F})$. The latter involution is denoted $\nu$. It has a simple geometric interpretation. Start in a point $x\in \,\,\Join
\! \!(\omega\vert_{\partial F})\subset \partial M$ and move along the   smooth immersed interval in $M$ adjacent to~$x$ and   formed by 
double points of $\omega(F)$. 
There are two possibilities for the second endpoint $y$ of this interval:  either  $y\in\, \,  \Join
\! \!(\omega\vert_{\partial F})$ and  then $ \nu(x)=y$ or $y\in B$ and then $\nu (x)=x$.

By the assumptions,  $\omega$ maps the   components  $K_1,\ldots , K_r$ of $\partial F$  to disjoint subsets of $\partial M$. Therefore \begin{equation}\label{i3} \Join
\! (\omega\vert_{\partial F})=\amalg_{i=1}^r \Join
\! (\omega (K_i)), \end{equation}  where $\Join
\! (\omega (K_i))$ is the set of all self-crossings of the loop $\omega(K_i)$.  Pick an arbitrary point $x\in \,\,
 \Join
\! (\omega\vert_{\partial F})$.  We   associate with $x$ a certain homology class as follows.  By \eqref{i3}, $x$  is   a self-crossing of   $\omega(K_i)$ for a unique  $i \in \{1,\ldots, r\}$.      Endow  $\omega(K_i)$ with  over/under-crossing  information  so that the resulting knot diagram on $\partial M$ presents $K_i\subset \partial M\times \RR$. Let $\omega^{-1}(x)=\{a,b\} $,  where  $a,b\in K_i $ and  we choose the notation so that
the  pair (the positive tangent vector  of $K_i$ at $a$, the positive tangent vector  of $K_i$ at $b$) is transformed by (the differential of) $\omega$ into a positive basis  in the tangent space of $\partial M$ at $x$.  Let   
$\gamma_x\subset K_i$ be the   positively oriented arc on $K_i$ leading from $a$ to $b$. The loop  $\omega_x=\omega(\gamma_x)$ is the \lq\lq distinguished half" of the diagram of $K_i$ at $x$.
Let 
 $[\omega_x]\in H_1(\partial M)$ be the homology class of $\omega_x$  and $\inc:H_1(\partial M)
\to H_1(M)$ be the inclusion homomorphism.   We now study the class $[\omega_x]$ in more detail.

Consider first the case where  $ \nu(x)=x$. We   claim    that then \begin{equation}\label{i1}\inc ([\omega_x]
)\in \omega_*(H_1(F))\subset H_1(M).\end{equation} Indeed,  the equality $ \nu(x)=x$ means that 
$\mu(a)\in
\{a,b\}$. Since $\mu$ is  fixed-point-free, 
$\mu(a)=b$.   By the definition of the involution $\tau:\widetilde \mmm \to \widetilde \mmm$, we have $\tau(a)=b$ and $
\tau(b)=a$.  Since       $\tau$ preserves the set of endpoints $\{a,  b\}$ of the interval $I_a\subset \widetilde \mmm$, we  necessarily  have  $\tau(I_a)=I_a$. The image of $I_a$ under the projection $p:\widetilde \mmm \to \mmm\subset F$ is an immersed interval on $F$ connecting $a$ and $b$. The product of the path  $\gamma_x \subset K_i\subset \partial F$   with the immersed  interval 
$p(I_a) $ oriented  from
$ b$ to $a$  is a loop in $F$, denoted  $\rho$. The loop  $\omega(\rho) $ in $M$  is a product of  
 the loops $\omega(\gamma_x)=\omega_x$  and  $\omega p\vert_{I_a}$. The latter  loop    has the form $ \delta 
\delta^{-1}$ where
$\delta $ is the path in $  M$ obtained by restricting
$ \omega\,p$ to the arc in $I_a$ leading from $b$ to   the unique   branch point of $\omega$  on $ I_a$.  Since the  loop  $ \delta 
\delta^{-1}$  is
contractible in $M$,  the loops $\omega_x$ and $\omega \rho $ are homologous in $M$. 
Since $\rho$ is a loop on $F$, we   conclude that 
$
\inc ([\omega_x]) \in \omega_*(H_1(F))$. 

Suppose now that $y=\nu(x)\neq x$. We shall establish two properties of the points $x,y\in \,  \Join
\! (\omega\vert_{\partial F})$. 
The first property is the   inclusion 
\begin{equation}\label{i2} \inc ([\omega_x]  + [\omega_{y}])\in  \omega_*(H_1(F) ) \subset H_1(M)  
\end{equation}
which we now verify.  As we know, $y\in \, \Join
\! (\omega (K_j))$ for some $j\in \{1,  \ldots, r\}$ (possibly $j=i$).
 Since the involution $\nu$ on $\partial \mmm/\tau=\partial \widetilde \mmm/\tau$ is induced by the involution  $\mu:  \partial \mmm\to    \partial \mmm$ and  $\omega^{-1}(x)=\{a,b\} $, we must have
 $\omega^{-1}(y)=\{\mu(a),\mu(b)\} $ with $\mu(a),\mu(b)\in K_j$. It is easy to check that  the  pair (the positive tangent vector  of $K_{j}$ at $\mu(a)$, the positive tangent vector  of $K_{j}$ at $\mu(b)$) is transformed by   $\omega$ into a {\it negative} basis  in the tangent space of $\partial M$ at $y$ (this change of orientation was  first pointed out in \cite{ca2}; cf.\ the argument in the next paragraph).
Therefore  the  oriented arc $\gamma_{y} \subset K_j$ begins  at    $\mu(b)$ and  
terminates  at $ \mu(a)$.
Consider   the  loop $\rho= 
\gamma_x \,p (I_{b})\,
\gamma_{y} \, (p(I_{ a}))^{-1}
$ in $F$ based at $a$. Here the intervals   $ I_{b}
,   I_{a}$ are   oriented  from
$b$ to $\mu(b)$ and from $a$ to $\mu(a)$, respectively.   Then   $\omega(\rho)$ is the product
of the loop $\omega_x$ based at $x$, the path $\omega p ( I_{b})$ beginning in $x$ and ending
in $ y$, the   loop $\omega_{y}$ based at $ y$, and  the path $(\omega p( I_{ a}))^{-1}$
beginning in
$ y$ and ending in $x$.   The  
paths 
$\omega p( I_{b})$ and   $(\omega p( I_{ a}))^{-1}$ are mutually inverse  since  $  I_{b} = \tau( I_{a}) $ and $\omega p \tau
=\omega p$. Hence 
 $ \inc ([\omega_x]  + [\omega_{y}]) =[\omega(\rho)]$. This implies   \eqref{i2}.

To state the second property of the pair $x,y=\nu(x)\neq x$, recall  the signs  of the self-crossings of a knot diagram introduced in    Section \ref{p2}. By definition,   $\sign(x)=+1$ if  the pair (the overgoing branch at $x$,  the undergoing branch at $x$) determines a positive basis in the tangent space of $\partial M$ at $x$ and $\sign(x)=-1$ otherwise.    We claim that 
  $\sign (x)=-\sign(y)$. Indeed, consider the smooth immersed  interval  $I=\omega p (I_a)$ in $M$  formed by 
double points of $\omega(F)$ and connecting  $x$ to $y$.  This interval is   an intersection of two branches of $\omega(F)$.  Each of these branches is  the projection of an immersed  ribbon  on $F$ whose core projects to $I$ and whose bases are short intervals on $\partial F$ parametrizing certain  branches of $\omega(\partial F)$ passing through $x$ and $y$.   
Since  the cores of the two ribbons in question project to the same interval $I$ in $M$, one of these   ribbons in $F\subset M\times \RR$ has to lie above the second one  with respect to the projection to $\RR$. This allows us to push the pair  (the overgoing branch at $x$,  the undergoing branch at $x$) along $I$ so that each branch is pushed along the corresponding ribbon transversally to the core. We  obtain at the end a pair (the overgoing branch at $y$,  the undergoing branch at $y$). If the first pair is positive with respect to the orientation of $\partial M$ induced from $M$, then the second one is negative and vice versa. Hence $\sign (x)=-\sign(y)$.

 We can now prove the lemma for $u_+$; the claim concerning $u_-$ can be proven similarly  or deduced by reversing orientation in $M$. Set $s_i=[\omega(K_i)]\in H_1(\partial M)$ for $i=1,\ldots ,r$ and
$s=  s_1+s_2+\cdots +s_r\in H_1(\partial M)$. By definition,  
\begin{equation}\label{sedrrttt}
\sum_{i=1}^r u_+ (K_i)= \sum_{i=1}^r   \,\,  \sum_{x \in
\Join  (\omega(K_i)), \, \omega_x\cdot  s_i   \neq 0, \,\sign ( \omega_x\cdot  s_i)=  \sign (x)}
\sign (x )\, 
t^{\vert  \omega_x\cdot s_i  \vert}.\end{equation}
   For  $x \in \,  \Join \!\!  (\omega(K_i))$,  the loop $\omega_x$ lies on
$\omega(K_i)$ and is      disjoint from $\cup_{j\neq i} \, \omega(K_j)$. Hence   
$ \omega_x\cdot  s = \omega_x\cdot  s_i$.  Using \eqref{i3}, we can rewrite  Formula  \ref{sedrrttt}  as
\begin{equation}\label{se968t}
\sum_{i=1}^r u_+ (K_i)= \sum_{x \in
\Join  (\omega\vert_{\partial F}), \,  \omega_x\cdot  s    \neq 0,\, \sign ( \omega_x\cdot  s )=  \sign (x)}
\sign (x )\, 
t^{\vert  \omega_x\cdot s  \vert}.\end{equation}

   We  now show  that each orbit of the involution $\nu$ on $\Join
\! \! (\omega\vert_{\partial F}) $ contributes $0$ to the right hand side of  \eqref{se968t}.  This will imply the claim of
the lemma. Consider first an orbit   consisting of one element $x\in\,  \Join
\! \! (\omega\vert_{\partial F})$ such that $\nu(x)=x$. 
It suffices to prove that    $ \omega_x\cdot  s =0$ because this would mean that $x$ does not contribute to the right hand side of   \eqref{se968t}.    Since the genus of $F$ is 0, the group $H_1(F)$ is generated by the homology classes of the boundary
components. Formula \eqref{i1} implies that there is  an integral  linear combination
$h $ of 
$s_1,\ldots,s_r
\in H_1(\partial M)$ such that  $\inc ([\omega_x])=\inc(h)$. Then 
$[\omega_x]- h \in K$ where $K$ is  the kernel of  the inclusion homomorphism $\inc: H_1(\partial M) \to H_1(M)$.  
The sum  $s=s_1+\cdots+s_r$, being represented by the boundary of   $ \omega(F)$,   also lies in $K$.  It is well known that  the homological intersection form on $H_1(\partial M)$ restricts to $0$ on $K$, that is
$ K\cdot  K =0$.  Hence 
$ ([\omega_x]- h )\cdot s=0$. Since the curves $\omega(K_1),\ldots,\omega(K_r)$ are pairwise disjoint and the intersection form is skew-symmetric,   $h \cdot s=0$. Therefore
$ [\omega_x] \cdot  s =0$.     

Consider   an orbit  of $\nu$  consisting of two distinct points $x, y\in\,  \Join
\! \! (\omega\vert_{\partial F})$. As we know, $\sign (x)=-\sign (y)$. 
Arguing  as in the previous paragraph with $[\omega_x]$ replaced by $[\omega_x] +[\omega_y]$ and using   \eqref{i2}, we conclude  that $ \omega_x\cdot  s =- \omega_{y} \cdot  s $. This   implies that either  
$ \omega_x\cdot  s =  \omega_{y} \cdot  s =0$ so that $x,y$ do not contribute to 
the right hand side of  \eqref{se968t} or the integers $ \omega_x\cdot  s $, $ \omega_{y} \cdot  s $ are non-zero and differ   by sign.  If $\sign ( \omega_x\cdot  s )=-\sign(x)$, then $\sign ( \omega_y\cdot  s )=-\sign(y)$ and $x,y$ do not contribute to the right hand side of  \eqref{se968t}.  If $\sign ( \omega_x\cdot  s )=\sign(x)$, then $\sign ( \omega_y\cdot  s )=\sign(y)$ and the contributions of $x,y$ cancel each other.  \end{proof}

\begin{corol}\label{aadddal:gg222} For any integers $m_1,\ldots,m_k\geq 1$, the polynomials
  $u_\pm^{(m_1,\ldots,m_k)}$ are cobordism invariants of knots.
		   \end{corol}

This follows   from  Lemma \ref{le} and Theorem \ref{th:0388}.

 \begin{corol}\label{aa+} For any slice knot $K$ and any integers $m_1,\ldots,m_k\geq 1$, we have   $u_\pm^{(m_1,\ldots,m_k)}(K)=0$. In particular,   $u_\pm (K)=0$.
		   \end{corol}
		   
		   For example, the knot    $K\subset S^1\times S^1 \times \RR$ constructed in Section \ref{exal}.2 is non-slice, since $u_+(K)=t\neq 0$.

  \section{Graded     matrices  and knots}\label{3}
                    
  We introduce  abstract  graded matrices and define  the  graded matrices of knots.

\subsection{Graded     matrices}\label{fi:g51}   A  {\it graded matrix} over 
an abelian group $A$  is a triple $(G, s, b  )$, where   $G$ is a finite set with distinguished element $s$ such that the set $G-\{s\}$ is partitioned as a union of two disjoint subsets $G_+$,  $G_-$ (possibly empty) and $b$ is an arbitrary mapping $G\times G \to A$. The partition $G-\{s\}=G_+\amalg G_-$ can be equivalently described in terms of the function  $\sign: G-\{s\}\to \{\pm 1\}$ sending $G_+$ to $+1$ and $G_-$ to $-1$.  An example of a   graded matrix 
is 
provided by the {\it trivial graded matrix} $(G,s, b)$ where $G=\{s\}$   and $b(s,s)=0$.

 Two graded matrices $(  G, s, b)$ and $(  G', s', b')$ over $A$ are {\it 
isomorphic} 
if there is a bijection $G\to G'$   
transforming $s$ into $s'$, $b$ 
into $b'$, and  $G_\pm$ into $G'_{\pm}$.

With a    graded matrix $T=(G,s,b)$ over   $A$, we associate  two graded matrices   $-T$ and $T^-$ over $A$. By definition,   $-T= (-G,s,-b)$, where $-G=G$ with    $(-G)_\pm=G_\mp$ and $(-b)(g,h)=-b(g,h)$ for  $g,h \in G$.   
 The graded matrix $T^-$ is defined by 
$T^-=(G,s, b^-)$,  where     $b^- (s,h)=-b( s,h)$, $b^- (h,s)=-b( h,s) $  for  $h\in G$, 
$b^-(g,h)= b(g,h)- b(g,s)-b(s,h) $ for  
$g,h\in G-\{s\}$, and the partition $G-\{s\}=G_+\amalg G_-$ is preserved.    The transformations  $T\mapsto -T$, $ T\mapsto T^-$ are commuting involutions on the class 
of graded matrices. 

\subsection{Skew-symmetric graded     matrices}\label{fi:g51+} 
A graded matrix $(G,s,b)$ over 
  $A$  is  {\it skew-symmetric} 
 if  the mapping
$b :G\times G  \to A$ is   skew-symmetric   in the sense that 
$b(g,h)=-b(h,g)$ for all $g,h\in G$  and  $b(g,g)=0$ for all 
$g\in 
G$.  

Given a skew-symmetric graded matrix $(G, s, b)$, we   say that an element $g\in G-\{s\}$ has {\it type 1}   if 
$b(g,h)=0$ for all $h\in G$. We say that $g\in G-\{s\}$ has {\it type 2} if 
$b(g,h)=b(s,h)$ for all $h\in G$.   We   call two elements $g_1,g_2\in 
G-\{s\}$ 
{\it complementary} if $\sign (g_1)=-\sign (g_2)$ and $  b (g_1, h)+ b
 ( g_2,h)=  b(s,h)$ for all $h\in  G$.  

We define  three 
transformations (or moves) $M_1, M_2, M_3$ on skew-symmetric graded matrices. The  moves $M_1, M_2, M_3$ delete   a type 1 element, a type 2 element, and  a pair  of complementary elements, 
respectively.  The inverse moves $M_1^{-1}, M_2^{-1}, M_3^{-1}$   add to a skew-symmetric  graded matrix $(  G, s, b)$ a type 1 element, a type 2 element, and  a pair  of complementary elements, 
respectively. 
More precisely, the move $M_1^{-1}$ (resp.\ $M_2^{-1}$)  transforms $(  G, s, b)$ into a  skew-symmetric graded matrix $(  
\overline G=G\amalg \{g\},  s, \overline b)$
such that    $\overline b:\overline G\times \overline G \to A$  extends
$b$ and  $\overline b (g,h)= 0$ (resp.\ $\overline b (g,h)=\overline  b (s,h)$) for all $h\in \overline G$.  The function $\sign$ on $\overline G -\{s\} $   extends the given function $\sign$ on $G-\{s\}$ and takes an arbitrary value $\pm 1$ on~$g$.
 The   move  $M_3^{-1}$  transforms $(  G, s, b)$ into a   graded matrix 
$(  \widehat 
G=G\amalg \{g_1,g_2\},  s, \widehat b)$
where   $\widehat b:\widehat G\times \widehat G \to A$ is any skew-symmetric   map   
extending 
$b$ and such that $\widehat b (g_1, h)+\widehat b
 ( g_2,h)=  \widehat  b(s,h)$ for all $h\in \widehat G$.  The function $\sign$ on $\widetilde G -\{s\} $   extends the given function $\sign$ on $G-\{s\}$,  takes an arbitrary value $\pm 1$ on $g_1$ and the opposite value on $g_2$.   
 Although we shall not need it, note that the   move  $M_2^{-1}$ can be expanded as a composition of $M_3^{-1}$ 
and  
$M_1$.

Two skew-symmetric graded matrices are {\it   homologous}
if they can be obtained from  each other by a finite sequence of 
  transformations $M^{\pm 1}_1,M^{\pm 1}_2, M^{\pm 1}_3$ and 
isomorphisms.   The homology is an equivalence relation on the class of skew-symmetric graded matrices over $A$.

\subsection{Primitive graded     matrices}\label{fi:g51+}    A graded matrix 
   over an abelian group  is  {\it primitive}    if it is skew-symmetric and has no   elements of type 1 or 2
  and 
no complementary pairs of elements.     For instance, the trivial graded matrix is primitive. 

Starting from an arbitrary skew-symmetric graded matrix $T  $ and 
recursively deleting    
  elements of types 1,  2  and   complementary pairs of 
elements,    we eventually obtain a  primitive graded matrix~$T_\bullet  $.  The following lemma shows that 
 $T_\bullet $
is determined by~$T$ uniquely up to  isomorphism.

\begin{lemma}\label{l:t51}   Two homologous primitive graded matrices are isomorphic.
		   \end{lemma}
                     \begin{proof}  We begin with the following assertion:
  
  $(\ast)$ a  move $M_i^{-1}$ followed by   $M_j$ with $i,j \in \{1,2, 3\}$ yields the same 
result as    
an isomorphism, or a   move $ M_k^{\pm 1}$, or     a  move $M_k $ 
followed 
by $M_l^{-1}$ with $k,l\in \{1,2, 3\}$.

  This assertion will imply the  lemma. Indeed,    
suppose that two primitive  graded matrices $T,T'$ are related by a finite 	
sequence 
of transformations $M_1^{\pm 1}, M_2^{\pm 1}, M_3^{\pm 1}$ and 
isomorphisms. 	
 An isomorphism   followed by $ M_i^{\pm 1}$     
can be also 
obtained as  $ M_i^{\pm 1}$    followed by an isomorphism.  Therefore 
all 
isomorphisms in our sequence   can be   accumulated at the end. The 
claim  
$(\ast)$ implies that $T,T'$
can be  related by a finite sequence of moves consisting of  several 
  moves of type $  M_i$ followed by  several 
  moves of type $ M_i^{-1}$ and   isomorphisms. However, since   $T$ is 
primitive  
  we cannot apply to it a move of type $  M_i $. Hence there are no 
such 
moves in the sequence. Similarly, since $
  T'$ (and any isomorphic graded matrix) is primitive, it cannot be obtained by 
an 
application of $M_i^{-1}$. Therefore our sequence  consists solely of 
isomorphisms so 
that $T$ is isomorphic to $T'$.

 Let us   prove $(\ast)$.    For $i, j \in \{1,2\}$, the move   $M_i^{-1}$ on a  skew-symmetric graded matrix 
$(G,s,b)$ adds 
one element $g$   and then  $M_j $ deletes   an  element  $g'\in G\amalg \{g\}$ of type $j$. If 
$g'=g$, 
then   $  M_j\circ M_i^{-1}=\id$. If   $g'\neq g$, then 
  the 
transformation $  
M_j\circ M_i^{-1}$   can be achieved by first applying $M_j $ that 
deletes 
$g'\in G$ and then applying $M_i^{-1}$ that adds $g$.

   Let $i\in \{1,2\}, j=3$. The move 
 $M_i^{-1}$ on $(G,s,b)$ adds an   element $g$ of type $i$ and  
  $M_j $ deletes   complementary  elements $g_1,g_2\in G\amalg \{g\}$. 
  If  $g_1,g_2\in G$, then $  
M_j \circ M_i^{-1}$ 
 can be achieved by first deleting $g_1,g_2$ and then adding $g$. If 
$g_1=g$, 
then   $g_2 \in G$ has type $3-i$  and $  M_j \circ M_i^{-1}$ is the 
move 
$M_{3-i} $    deleting $g_2$. The case $g_2=g$ is similar.

   Let $i=3, j\in \{1,2\}$. The move 
 $M_i^{-1}$  on $(G,s,b)$ adds   complementary  elements $g_1,g_2$   and  
$M_j  $ deletes a certain $g\in G\amalg \{g_1,g_2\}$ of type $j$. 
  If   $g\in G$, then  $  M_j \circ 
M_i^{-1}$  can 
be achieved by first deleting $g$ and then adding $g_1,g_2$.
   If $g=g_1$, then   $ g_2$ has type $3-j$ 
and $  
M_j \circ M_i^{-1}=M_{3-j}^{-1}$. The case $g=g_2$ is similar.

      Let $i=j=3$. The move 
 $M_i^{-1}$  on $(G,s,b)$ adds   complementary  elements $g_1,g_2$   and  
$M_j $ deletes    complementary elements  $g'_1,g'_2\in G\amalg \{g_1,g_2\}$. If the pairs  $g_1,g_2$   and  $g'_1,g'_2$   are disjoint, then 
  $  M_j \circ M_i^{-1}$  can be achieved by first deleting  $g'_1,g'_2\in G$ 
and then 
adding $g_1,g_2$. If  these 
two pairs     
  coincide, then $  M_j \circ M_i^{-1}$  is the identity. It remains to 
consider 
the case where these pairs have one common element, say $g'_1=g_1$, 
while 
$g'_2\neq g_2$. Then $g'_2\in G$ and for all $h\in G$,
  $$\widehat b(g_2,h)= \widehat  b(s,h)-\widehat b(g_1,h)= \widehat  b(s,h)-\widehat b(g'_1,h)=\widehat 
b(g'_2,h)= 
b(g'_2,h).$$
 Therefore the move   $  M_j \circ M_i^{-1}$ produces a graded matrix isomorphic 
to 
$(G,s,b)$. The isomorphism $  (G-\{g'_2\}) \cup \{g_2\} \approx G$ is the 
identity on 
$G- \{g'_2\}$ and  sends $g_2$ into $g'_2$. Note that $\sign (g'_2)=-\sign (g'_1)=-\sign (g_1) =\sign (g_2)$.
\end{proof}

Note finally that the transformations  $T\mapsto -T$, $ T\mapsto T^-$    preserve the class of   skew-symmetric graded matrices and      are   compatible
with the relation of  homology on this class. Clearly,     $(-T)_\bullet =- T_\bullet$ and $(T^-)_\bullet=(T_\bullet)^-$.

 \subsection{Graded matrices of knots}\label{fi:g53} With a  
knot diagram $ D$ on a surface~$\Sigma$ we associate a skew-symmetric  graded matrix 
$T( D)=(G,s,b)$ over $\ZZ$ as follows. Set $G =\{s\} \,\amalg 
\Join \!\! (D)$ and provide   $G-\{s\}=\,\Join \!\! (D)$ with the bipartition  $G_+=\,\Join_+ \!\! (D)$ and $ G_-=\,\Join_- \!\! (D)$.  To define the pairing
  $b :G\times G\to \ZZ$, consider the mapping  $\alpha: G\to H_1(\Sigma)$ sending $s$ to $[D]\in H_1(\Sigma)$ and sending any $x\in \, \Join \!\! (D)$ to   $[D_x]\in H_1(\Sigma)$.    The pairing
  $b $  is defined by $b(g,h)=\alpha(g) \cdot \alpha (h) \in \ZZ$ for  $g,h \in G$.  It is clear that $b$ is skew-symmetric.

 \begin{lemma}\label{th:e53}  If two   knot diagrams   on $\Sigma$ present isotopic knots in $\Sigma \times \RR$, 
then their 
graded matrices   are homologous. \end{lemma}
  \begin{proof}  It suffices to  verify that, if   knot diagrams  $D, D'$ on $\Sigma$ are related by Reidemeister moves, 
then their 
graded matrices   are homologous. The arguments given in Section \ref{1.3} in the proof of isotopy invariance of $u_\pm$ shows that   if $D'$ is obtained from $D$ by   
the first (resp.\ the second) R-move, then $T(D)$ is obtained from  
$T(D')$ by $M_1$ or $M_2$ (resp.\ by $M_3$).  If $D'$ is obtained from $D$ by   
the third  R-move, then $T(D)$ is  isomorphic to $T(D')$.   \end{proof} 
  
This lemma allows us to define the {\it graded matrix of a    knot} $K\subset \Sigma \times \RR$. Present $K$ by a  diagram $D$ on $\Sigma$. By Lemma \ref{th:e53}, the homology class of the skew-symmetric graded matrix $T(D)$ does not depend on the choice of $D$. By Lemma \ref{l:t51},  
the  
primitive graded matrix  over $\ZZ$ defined  by $T_\bullet(K)=(T(D))_\bullet$ is determined by $K$ uniquely up to isomorphism.  The isomorphism class of  $T_\bullet(K)=(G_\bullet, s_\bullet, b_\bullet) $  is a diffeomorphism invariant of $K$.      Note one application of this  invariant:   
  any diagram of $K$ on  $\Sigma  $   must have at least $\card (G_\bullet ) -1 $ double points. A related fact:  if $K$ is presented by a diagram $D$ such that  $T(D)$ is primitive, then 
  any diagram of $K$ on $\Sigma$ has at least as many double points as $D$. Indeed, in this case $T_\bullet (K)=T(D)$ and $\card (G_\bullet ) -1=\card (\Join \!\! (D))$.

         It is easy to check that       $T_\bullet(-K)= (T_\bullet(K))^-$  and    $T_\bullet(\overline K)= -(T_\bullet(K))^-$.   If    $K$ bounds a disk in $\Sigma\times \RR$, 
then 
$T_\bullet(K)$ is the trivial graded matrix.  

More generally, for a    finite sequence of  positive integers $m_1,\ldots, m_k  $, we   define the     {\it higher  graded matrix}  of $K$ by  
$$T_\bullet^{m_1,\ldots, m_k }  ( K)= T_\bullet(  K^{(m_1,\ldots,m_k)}) .$$
This primitive graded matrix over $\ZZ$   is a diffeomorphism invariant of $K$.

\subsection{Polynomials $u_\pm$ re-examined}\label{oper}  For a    graded matrix
$T=(G,s,b)$ over $\ZZ $, set    $$u_\pm (T)=\sum_{x\in G, \, b(g,s)  \neq 0,\, \sign (b(g,s))=\pm \sign (x) } \sign (x) \,t^{ \vert  b(g,s) \vert }\in \ZZ[t] \,.$$
It is clear that $u_\pm (-T)=-u_\pm (T)$ and $u_\pm (T^-)= u_\mp (T)$.  
Both $u_+ $ and $u_- $ are    homology invariants of
skew-symmetric graded matrices. 
 
The polynomial invariants $u_\pm (K)$ of any knot $K$ can be computed  from  $T_\bullet(K)$ by    $u_\pm(K)= u_\pm(T_\bullet(K))$.  Indeed,  for a   diagram $D$ of $K$,
  $$  u_\pm(K)=u_\pm(T(D))=u_\pm((T(D))_\bullet)=u_\pm(T_\bullet(K)) .$$ 
  More generally,  
  $u_\pm^{m_1,\ldots, m_k }  ( K)=u_\pm (T_\bullet^{m_1,\ldots, m_k }  ( K))$ for any integers $m_1,\ldots, m_k \geq 1$.
  
   \subsection{Remark}\label{fierf55}   By \cite{tu2}, any loop on a surface gives rise to a   based matrix (based matrices are defined as graded matrices $(G,s,b)$   but   no  bipartition of $G-\{s\}  $  is distinguished). The based matrix of the loop  underlying a knot diagram $D$   is obtained from     $T(D) $ by  forgetting the  bipartition.

  \subsection{Examples}\label{exal+} 1.  Consider the knot $K$ on $  S^1\times S^1 $ presented by  the  diagram $D$   with two crossings $x,y$ from  Example \ref{exal}.2.  Then $T(D)=(G, s, b)$, where $G=\{s, x,y\}$ with $G_+=\{x,y\} , G_-=\emptyset$ and  the mapping $b:G\times G \to \ZZ$ is given by the   matrix
\begin{equation}\label{mal}
 \left [  \begin{array}{ccccc}  0& -1&1 \\
        1& 0& 1  \\
-1& -1& 0 
\end{array} \right ],\end{equation}
where the rows and the columns correspond to $s,x,y$, respectively. This matrix is primitive and therefore $T_\bullet (K)=(T(D))_\bullet=T(D)$. The equalities $D_x\cdot D=-D_y\cdot D=1$ imply that the   knot  $K^{(m)}$ with $m\geq 2$  is presented by a diagram  without self-crossings;  therefore the higher graded matrices of $K$ are trivial.   

2.  For any   integers $p,q\geq 1$, the author defined in \cite[Sections 3.3 and  4.3]{tu2} a closed curve $\alpha_{p,q}$ on a closed (oriented) surface $\Sigma$. This curve is defined by its Gauss diagram   consisting of a     circle in $\RR^2$ with $p$ horizontal chords directed leftward and $q$ vertical chords directed upward  (each   horizontal chord  should cross each   vertical chord    inside the circle). The genus of $\Sigma$   is  equal to  1 
if 
$p=q=1$, to  3 if $\min(p,q)\geq 3$, and to $2$ in all the other cases. The   curve $\alpha_{p,q}$ has $p+q$ double points $x_1,\ldots, x_{p+q}$ corresponding to the   chords of the Gauss diagram. To transform $\alpha_{p,q}$ into a  diagram of a knot $K$ on $\Sigma$, it is enough to fix a  function $\sign: \{x_1,\ldots, x_{p+q}\}\to \{\pm 1\}$.  
 The corresponding skew-symmetric graded matrix $T=(G,s,b:G\times G \to \ZZ)$ is  computed as follows (cf. \cite{tu2}): $G=\{s, x_1,\ldots, x_{p+q}\}$; 
      $b(x_i,s) =q$ for $i\in \{1,\ldots, p\}$ and $b(x_{p+j},s)= -p$ for 
$j\in \{1,\ldots, q\}$; $b(x_i,
x_{i'}) = b(x_{p+j},
x_{p+j'}) =0$ for   $i,i'\in \{1,\ldots, p\}$, $j,j'\in \{1,\ldots, q\}$;
and  $b(x_i,
x_{p+j})= p+q+1-i -j  $ for $i\in \{1,\ldots, p\}$, 
$j\in \{1,\ldots, q\}$. The partition $G-\{s\}=G_+\cup G_-$ is determined by the function $\sign$ in the usual way. It is easy to check from the definitions that the graded matrix $T$ is primitive except in the case where $p=q=1$ and $\sign(x_1)=-\sign (x_2)$. Excluding this case, we obtain that  $T_\bullet (K)=T$. This implies in particular that $K$ is non-trivial and any   diagram on $\Sigma$ presenting a knot isotopic to $K$ has at least  $p+q$ double points. To compute the polynomial $u_\pm(K)$, set $$a_\pm= \card \{1\leq i \leq p, \sign (x_i)=\pm 1 \},$$
$$b_\pm= \card \{1\leq j \leq  q, \sign (x_{p+j})=\pm 1 \}.$$ Then  
$$u_\pm (K) =u_\pm(T)=\pm (a_\pm \,t^{q } -  b_\mp  \,t^{p }) .$$
  This formula easily implies  that such knots  $K$ corresponding to different pairs $(p,q)$ are never diffeomorphic and that two knots corresponding to the same pair $(p,q)$ and different functions  
  $ \{x_1,\ldots, x_{p+q}\}\to \{\pm 1\}$ may be isotopic only if these two  functions take the value $+1$ the same number of times and the restrictions of 
these   functions to the set $\{x_1,\ldots, x_p\}$ take the value $+1$ the same number of times. 
In the case  where $p=q=1$ and $\sign(x_1)= \sign (x_2)=1$ we recover the knot $K$ from the previous example.  In the exceptional case $p=q=1$ and $\sign(x_1)=-\sign (x_2)$ the knot $K$ is isotopic to a knot on $\Sigma$ presented by a diagram without self-crossings.  

In this example, the higher graded matrix $T^m_\bullet (K)$ is non-trivial if the integer $m $ divides $p$ or $q$. For instance, if $m$ divides both $p$ and $q$, then   $K^{(m)}=K$ and $T^m_\bullet (K)=T$. If $m$ divides $p$ but not $q$, then $T^m_\bullet (K)$ is  obtained from $T$ by deleting $x_1,\ldots, x_p$.

     \section{Genus and cobordism for  graded matrices}\label{5}

We introduce a numerical genus for finite families of  graded matrices  and use it to define    a relation of cobordism for graded  matrices (not necessarily skew-symmetric). This relation will be confronted with knot  cobordism      in the next section.

Throughout this section the symbol $R$ denotes a domain, i.e., a commutative ring without 
zero-divisors. By a  graded matrix 
 over $R$, we   mean a graded matrix over the additive group of $R$.

\subsection{The genus}\label{uussstrg52}   Consider  a family of $r\geq 1$ graded matrices 
$T_1=(G_1,s_1,b_1)$, $\ldots$, $T_r=(G_r,s_r,b_r)$ over $R$.  
We define a numerical invariant $\sigma (T_1,\ldots, T_r)$ of this family called the {\it graded genus}.    First of all, replacing $T_1,\ldots,T_r$ by isomorphic graded
matrices, we can assume that the sets
$G_1,\ldots,G_r$ are disjoint.  
Set $G=\cup_{t=1}^r G_t$ and $G_\pm= \cup_{t=1}^r (G_t)_{\pm}\subset G$. Denote $RG$ the free $R$-module  with basis  $G$. The maps $\{b_t:G_t\times G_t\to R\}_t$ induce a  bilinear form $b=\oplus_t b_t:RG\times
RG\to
R$ such that 
$b(g,h)=b_t(g,h)$ for 
$g,h\in G_t$  and 
$b(G_t,G_{t'})=0$ for $t\neq t'$.

Let $S$ be the submodule of  $RG$ generated by $s_1,\ldots,s_r$.   We call a vector
$x\in
RG$ {\it short} if   $x\in S$  or $x\in g+
S$ for some  $g\in 
G-\{s_1,\ldots,s_r\} $  or $x \in g+h+
S$ for      $g,h \in G-\{s_1,\ldots,s_r\}$ of opposite sign (that is  $g\in G_+, h\in G_-$ or $g\in G_-, h\in G_+$).   
  A  {\it (graded) filling} of
$T_1,\ldots,T_r$ is  a finite family $   \{\lambda_i\}_i$ of short vectors in $RG$ such that $\sum_i
\lambda_i=\sum_{g\in G} g\,  (\modu  S)$ and one of $\lambda_i$ is equal to $s_1+s_2+\cdots+s_r$.
Note that each basis vector $g\in G-\{s_1,\ldots,s_r\}$ appears in exactly one $\lambda_i $ with  coefficient 
$+1$ and does not appear in   other $\lambda_i$'s. The basis vectors
$ s_1,\ldots,s_r $  may appear in several $\lambda_i$ with non-zero coefficients.  
For example, the set consisting of all elements of $G $ and the vector  $s_1+s_2+\ldots+s_r$
is a filling of $T_1,\ldots,T_r$.

The {\it matrix}
of  a filling 
$\lambda = \{\lambda_i\}_i$    is the square matrix
$(b(\lambda_i,\lambda_j))_{i,j}$ over
$R$. Let  $\sigma(\lambda)  $ be half of its  rank.  This number  is equal to   half of the  rank of the restriction of $b$  to the submodule   of $RS$ generated by the vectors $ \{\lambda_i\}_i$. The number $\sigma(\lambda)\geq 0$ is an integer or a half-integer;  it is certainly an integer if 
$T_1,\ldots,T_r$ (and then $b$) are skew-symmetric. Set  $$\sigma(T_1,\ldots,T_r)=\min_\lambda
\sigma(\lambda)\geq 0, $$ where $\lambda$ runs over all fillings of $T_1,\ldots,T_r$.  
Clearly   $\sigma (T_1,\ldots,T_r)=0$ if and only if $(T_1,\ldots,T_r)$ has a filling with
zero matrix.  In this case we call the family  
$T_1,\ldots,T_r$  {\it hyperbolic}.

It is obvious that  
$\sigma(T_1,\ldots,T_r)$  is preserved when   $T_1,\ldots,T_r$ are permuted or 
replaced with isomorphic graded matrices.   It is easy to check that 
$$\sigma(-T_1,\ldots,-T_r)=  \sigma(T_1,\ldots,T_r).$$  If
$T_r$ is a trivial graded matrix, then
$\sigma(T_1,\ldots, T_r)=\sigma(T_1,\ldots,
T_{r-1})$ (because then the vector $s_r\in S$ lies in the annihilator of the form $b$).

For $r=1$, the definitions above apply to a single graded matrix $T$ and yield  the genus $\sigma(T)\geq 0$.  By definition, $T$ is hyperbolic if and only if $\sigma (T)=0$.
  
   We state   a key property of the genus.

\begin{lemma} \label{ccxf} For any $1\leq t\leq r$ and any graded matrices  $T_0, T_1,\ldots, T_r$,
$$\sigma(T_1,\ldots, T_r)\leq  \sigma(T_1,\ldots, T_t, T_0)
+ \sigma(-T_0, T_{t+1},\ldots, T_r).$$
\end{lemma}
\begin{proof} Consider for concreteness the case where $t=1$ and $r=2$, the general case is   similar.  
We must prove that $\sigma(T_1, T_2)\leq  \sigma(T_1,  T_0)
+ \sigma(-T_0, T_{2})$.  Let
$T_i=(G_i,s_i,b_i)$ for $i=0,1,2$ and   $T'_0=(G'_0,s'_0,b'_0)$ be a  copy of $T_0$ where 
$G'_0=\{g' \,\vert \, g\in G_0\}$, $s'_0=(s_0)'$,
and $b'_0$ is defined by $b'_0(g',h')=b_0(g,h)$ for $g,h\in G_0$.   We can assume
that the sets
$G_1, G_0,G'_0, G_2$ are disjoint.  Let  $\Lambda_1 , \Lambda_0, \Lambda'_0, \Lambda_2 $ be  free $R$-modules
freely generated by
$G_1, G_0, G'_0,  G_2$, respectively, and $\Lambda=  \Lambda_1\oplus  \Lambda_0\oplus \Lambda'_0\oplus 
\Lambda_2 $.
There is  a  unique skew-symmetric bilinear form $B=b_1\oplus b_0\oplus (-b'_0)\oplus b_2$ on $\Lambda$ such that the
sets
$G_1, G_0, G'_0, G_2\subset \Lambda$ are mutually orthogonal and the restrictions of $B$ to these  sets  are equal
 to $b_1, b_0, -b'_0,b_2$, respectively.

 Let $\Phi$   be the submodule of
$ \Lambda_0\oplus \Lambda'_0$ generated by    the vectors
$\{g+g'\}_{g\in G_0}$.  Set $L= \Lambda_1\oplus  \Phi \oplus \Lambda_2\subset \Lambda$.     The
projection $p:
L  \to \Lambda_1\oplus    \Lambda_2$  along $\Phi$ transforms $B$ into $b_1\oplus b_2$. 
Indeed, for any  $g_1,h_1\in G_1, g , h\in G_0, g_2,h_2\in G_2$, 
$$B(g_1+g  +g' +g_2, \,\,h_1+ h + h'+h_2)$$
$$=
b_1(g_1 ,h_1)+ b_0 (g ,  h)+ (-b'_0) (g', h')+ b_2(g_2,h_2)= 
b_1(g_1 ,h_1)+  b_2(g_2,h_2).$$

Pick    a filling  $ \{\lambda_i\}_i\subset \Lambda_1 \oplus \Lambda_0$ of the pair 
$(T_1, T_0) $ whose matrix has rank   $2\sigma(T_1,  T_0)$.  This means that   the restriction of $B$ to the
 module $V_1\subset \Lambda_1\oplus \Lambda_0$ generated by $ \{\lambda_i\}_i$ has rank   $2\sigma(T_1,  T_0)$. 
Similarly, pick      a filling
$ 
\{\varphi_j\}_j\subset
\Lambda'_0\oplus \Lambda_2$    of the pair  
$(-T'_0, T_2) $ such that 
  the restriction of $B$ to the
 module $V_2\subset \Lambda'_0\oplus \Lambda_2$ generated by $ \{\varphi_j\}_j$ has rank  $2\sigma (-T'_0, T_2)$.
  We claim that there is a  finite  set $\psi\subset (V_1+V_2)\cap L$ such that  the set $p(\psi) \subset
\Lambda_1\oplus    \Lambda_2$ is a filling of    the pair  
$(T_1, T_2) $. Denoting by
$V$ the submodule of   
$\Lambda_1\oplus    \Lambda_2$ generated by $p(\psi)$, we obtain  then   the desired inequality:
$$\sigma(T_1, T_2)\leq  \sigma (p(\psi))= (1/2)  \rk ((b_1\oplus b_2) \vert_V)= (1/2)  \rk (B \vert_{p^{-1}(V)})
$$
$$ \leq 
(1/2) 
\rk (B
\vert_{(V_1+V_2)\cap L} )
\leq   (1/2) \rk (B
\vert_{V_1+V_2})$$
$$=  (1/2)  \rk (B \vert_{V_1}) + (1/2) \rk (B \vert_{V_2})=\sigma(T_1,  T_0)
+ \sigma(-T_0, T_{2}).$$
Here the second inequality follows from the inclusion $$p^{-1}(V)\subset  (V_1+V_2)\cap L + \Ker p$$ and the fact
that
$\Ker p=\Phi$ lies in the annihilator of $B\vert_L$.

To construct $   \psi$, we first modify     $ \{\lambda_i\}_i$  as
follows.  Let $\lambda_1$ be the vector  of this   filling equal to $s_1+s_0$.  Adding appropriate 
multiples of
$\lambda_1$  to other  
$\lambda_i$ we can ensure that the basis vector $s_0 \in G_0 $ appears in all $\{\lambda_i\}_{i\neq
1}$  with coefficient 0.  This transforms  $ \{\lambda_i\}_i$ into a new filling of  $(T_1, T_0) $ which will be
from now on denoted    $\lambda= \{\lambda_i\}_i$. This transformation does not change the module
$V_1$ generated by the vectors of the filling.   Similarly, we can assume that  a vector  $\varphi_1$ of the filling
$\varphi=\{\varphi_j\}_{j }$ is equal to 
$s'_0+s_2$ and  that the basis vector $s'_0 \in G'_0 $ appears in all
$\{\varphi_j\}_{j\neq
1}$ with coefficient~0.

The
filling $\lambda$ gives rise to an oriented  1-dimensional manifold     $\Gamma_\lambda$ with boundary $(G_1-\{s_1\})\cup  (G_0-\{s_0\})$.   It is defined as follows.      Each
$\lambda_i $ having the form 
$g+h \, (\modu R s_1)$  with $g,h\in  (G_1-\{s_1\})\cup  (G_0-\{s_0\})$ gives rise to a component of
$\Gamma_\lambda$ diffeomorphic to $[0,1]$ and connecting   $g$ with $h$. By the definition of a filling,   the elements $g,h$ have opposite signs. We orient the component   in question so that it leads from the element with sign $-1$ to  the element with sign $+1$.  Each $\lambda_i $ having
the form 
$g  \, (\modu R s_1 )$ with $g \in  (G_1-\{s_1\})\cup  (G_0-\{s_0\})$ gives rise to a component of
$\Gamma_\lambda$ which is a copy of $[0, \infty)$  
where $0$ is identified with  $g$.  We orient this component towards $g$ if   $\sign (g)=+1$ and out of $g$ if $\sign (g)=-1$. All the other $\lambda_i$ and in particular $\lambda_1$ do not contribute to
$\Gamma_\lambda$.  The definition of a filling implies that
$\partial
\Gamma_\lambda= (G_1-\{s_1\})\cup  (G_0-\{s_0\})$.   Similarly,   the filling
$\varphi$ gives rise to an oriented     1-dimensional manifold 
$\Gamma_\varphi$ with boundary 
$(G_0-\{s'_0\})\cup  (G_2-\{s_2\})$. We can assume that  $\Gamma_\lambda$ and 
$\Gamma_\varphi$ are disjoint.   Gluing  $\Gamma_\lambda$ to
$\Gamma_\varphi$ along the canonical identification
$G_0-\{s_0\} \to G'_0-\{s'_0\},\, g\mapsto g'$, we obtain a 1-dimensional manifold  $\Gamma$  with $\partial
\Gamma=(G_1-\{s_1\})\cup  (G_2-\{s_2\})$.  Note that the orientations of  $\Gamma_\lambda$, 
$\Gamma_\varphi$ at $G_0-\{s_0\} = G'_0-\{s'_0\}$ are compatible because $(G'_0)_\pm =(G_0)_\mp$ (cf.\ the definition of   the transformation $T\mapsto -T$). Therefore these orientations extend to an orientation of $\Gamma$.   Each component  $C$ of $\Gamma$ is glued from several components of 
$\Gamma_\lambda\amalg  \Gamma_\varphi$ associated with certain vectors  $$\lambda_i\in V_1\subset  \Lambda_1\oplus
\Lambda_0
\subset \Lambda\quad {\text {or}} \quad \varphi_j
\in V_2\subset \Lambda'_0\oplus \Lambda_2 \subset \Lambda.$$ Let   $\psi_C\in V_1+V_2 \subset 
\Lambda$ be the sum of these vectors over all components of 
$\Gamma_\lambda\amalg  \Gamma_\varphi$ contained in $C$.  The following two facts imply that  $\psi_C\in
L$: 

 (i) each point of $C\cap (G_0-\{s_0\})\approx C\cap (G'_0-\{s'_0\})$ is adjacent to   one component of  
$\Gamma_\lambda$ and 
 to one component of   $\Gamma_\varphi$ and 
 
 (ii) $s_0  $ does not show up   in  
$\{\lambda_i\}_{i\neq 1}$  and   $s'_0   $ does not show up   in  
$\{\varphi_j\}_{j\neq
1}$. 

Set
$\psi_1=\lambda_1+\varphi_1=s_1+s_0+s'_0+ s_2
\in
\Lambda$. Clearly, $\psi_1 \in  (V_1+V_2)\cap L$.    
Set
$\psi=\{\psi_1\}\cup \{\psi_C\}_C$ where $C$ runs over the  components of $\Gamma$ with non-void boundary. As we know, $\psi\subset (V_1\oplus V_2) \cap C$. 
 Let us check that $p(\psi) \subset
\Lambda_1\oplus    \Lambda_2$ is a filling of   
$(T_1, T_2) $.  Observe that for a compact component $C$ of $\Gamma$ with endpoints $g,h\in G_1\cup G_2$, we have 
$p(\psi_C)=g+h \, (\modu R s_1+R s_2)$.  If $C$ is oriented, say, from $g$ to $h$, then necessarily $\sign (g)=-1$ and   $\sign (h)=+1$. For a non-compact component $C$ of $\Gamma$ with  one endpoint 
$g \in G_1\cup G_2$, we have 
$p(\psi_C)=g \, (\modu R s_1+R s_2)$. In both cases the vector $p(\psi_C)$ is short. The sum of all  vectors in the family $\psi $   is equal to
$\sum_{g\in G_1\cup G_2} g \, (\modu R s_1+R s_2)$.  Also $p(\psi_1)=s_1+s_2$.  This means that $p(\psi)$ is a
filling of 
$(T_1, T_2)
$ so that $\psi$ satisfies all the required conditions.
\end{proof}

\subsection{Cobordism of graded matrices}\label{u125452} Two graded matrices $T_1,T_2$ over $R$  are {\it cobordant} if
$\sigma(T_1,-T_2)=0$.

\begin{theor} \label{859s}  (i) Cobordism is an equivalence relation on the   class of graded matrices. 

(ii) Isomorphic graded matrices are cobordant.

(iii)  The graded genus   of a family of graded matrices  is a cobordism invariant.

(iv)  A graded matrix   is cobordant to the  trivial graded matrix if and only if  it  is hyperbolic. 

(v) Homologous skew-symmetric graded matrices are cobordant.
\end{theor}
\begin{proof} (i) and (ii)   For a  graded matrix $T=(G,s,b)$, the graded matrix $-T$ is isomorphic to the triple $ (G',s',b')$  where
$G'=\{g'\, \vert \, g\in G\}$ is a  disjoint copy of $G$, $(G')_\pm=G_\mp$,  and  $b'(g',h')=-b(g,h)$ for any $g,h\in G$.  The vectors $\{g+g'\}_{g\in G}$ form a  filling
 of  the pair $(T, -T)$.  The matrix of this  filling is
0. Therefore 
$\sigma(T,-T)=0$ so that $T$ is cobordant to itself.  A similar argument proves (ii). That the relation of cobordism is symmetric follows from the equalities
$\sigma(T_2,-T_1)=\sigma (-T_2,T_1)=\sigma(T_1,-T_2)$.  The transitivity    follows from the inequalities
$$0\leq \sigma (T_1,-T_3)\leq  \sigma (T_1,-T_2)+\sigma (T_2,-T_3)$$
which is a  special case of Lemma \ref{ccxf}.

(iii) We need to prove that    $\sigma(T_1,\ldots,T_r)$ is preserved  when  $T_1,\ldots,T_r$ are replaced with
cobordant graded matrices.  By induction, it suffices to prove that 
\begin{equation}\label{fuy}\sigma(T_1,\ldots,T_{r-1} , T'_r)=\sigma(T_1,\ldots,T_{r-1} , T_r)\end{equation} for any graded matrix   $T'_r$  cobordant to $T_r$.   
By Lemma \ref{ccxf}, 
$$\sigma(T_1,\ldots,T_{r-1} , T_r)\leq \sigma(T_1,\ldots,T_{r-1} , T'_r)+\sigma(-T'_r ,
T_r) =\sigma(T_1,\ldots,T_{r-1} , T'_r).$$  Similarly, 
$\sigma(T_1,\ldots,T_{r-1} , T'_r)\leq 
\sigma(T_1,\ldots,T_{r-1} , T_r)$. This implies \eqref{fuy}. 

(iv)  If  a  graded matrix $T$ is cobordant to the  trivial graded matrix $T_0=(\{s_0\}, s_0$, $b=0)$, then $\sigma
(T)=\sigma(T_0)=0$ and  
$T$ is hyperbolic. Conversely, if $T=(G,s,b)$ is hyperbolic, then it has a filling $\lambda=\{\lambda_i\}_i$ with zero matrix. Since $s=\lambda_i$ for some $i$, we necessarily have $b(s, \lambda_i)=0$ for all $i$ (and in particular, $b(s,s)=0$). Adding to $\lambda$  the vector $s+s_0$, we obtain a filling of the pair $(T, T_0)$ with zero matrix. Hence $T$ is cobordant
to $-T_0=T_0$.

(v) Let a skew-symmetric graded matrix $T' =(G' ,s' ,b' )$  be  obtained from a skew-symmetric graded matrix $T=(G,s,b) $ by a move $M_i^{-1}$
with
$i=1,2,3$. We can assume that   $G'$  is a  union of a disjoint copy $ \{h'\,\vert \,h\in G\}$ of
$G$ and one new element $g$  in the case $i=1,2$ or two new elements $g_1,g_2$ in the case 
$i=3$.   For
$i=1$ (resp.\
$i=2$, $3$),  the vectors
$\{h+h'\}_{ h\in G}$  and the vector $g$ (resp.\ $g-s'$, $g_1+g_2-s'$)  form a filling of  the pair $(T,
-T')$. The matrix of this   filling  is zero. Hence $\sigma(T,-T')=0$ so that $T$ is cobordant to
$T'$. 
\end{proof}

\begin{corol}\label{aaal:gfglpmns5} For any  skew-symmetric  graded matrices  $T_1,\ldots,T_r$,  
 $$\sigma((T_1)_\bullet,\ldots,(T_r)_\bullet)=\sigma(T_1,\ldots, T_r).$$
		   \end{corol}

 \subsection{Remarks}\label{sdm1}   We point out a few  further  properties of   graded matrices.
 
 (a)  Let us call a filling  of   a graded matrix $T=(G,s,b)$  {\it simple} if all vectors of the filling are pairwise distinct and  either belong to the set of basis vectors $G$ or have the form   $g+h$ with $g, h\in G-\{s\}$ and $\sign (g)=-\sign (h)$. Any filling $\lambda$ of $T$ can be transformed into a simple one by adding     vectors proportional to $s$ to the vectors of $\lambda$ distinct from~$s$ (and eliminating repetitions).  This transformation does not change the module generated by the vectors of $\lambda$ and does not change   $\sigma(\lambda)$. Thus,  to compute $\sigma (T)$ we can    restrict ourselves to the simple fillings. The simple fillings of $T=(G,s,b)$ are finite in number; they  bijectively correspond to involutions $\nu$ on $G-\{s\}$ such that
   every free orbit of $\nu$ meets both $G_+$ and $G_-$. In particular, if one of the sets $G_+, G_-$ is empty, then $T$ has only one simple filling consisting of all elements of $G$ and then $\sigma (T)$ is half of the rank of the matrix $(b(g,h))_{g,h\in G}$. Such a graded matrix $T$ is hyperbolic if and only if $b=0$. For example,  the graded matrix of Section \ref{exal+}.1   satisfies $G_-=\emptyset$, $b\neq 0$ and therefore it  is not hyperbolic.

 (b)  The function  $(T_1,T_2)\mapsto \sigma(T_1,-T_2)$ defines a metric on the set of cobordism classes of graded 
matrices over a domain.

(c) For  any graded matrices  $  T_1,\ldots, T_r$, $T'_1,\ldots,T'_q$ with $q\geq 1$ and any  $ 1\leq t\leq r$,  
$$\sigma(T_1,\ldots, T_r)\leq  \sigma(T_1,\ldots, T_t, T'_1,\ldots,T'_q)
+ \sigma(-T'_1,\ldots,-T'_q, T_{t+1},\ldots, T_r)+q-1.$$

(d)  For a graded matrix $T=(G,s,b)$ over a domain $R$, the element $b(s,s)$ of~$R$  is a cobordism invariant of $T$. We call  $T $ {\it normal} if $b(s,s)=0$. For example, all  skew-symmetric graded matrices   are   normal.  For normal graded matrices over $\ZZ$,  the polynomials $u_+$ and $u_-$   are cobordism invariants.  This follows from the following more general fact: if a family $T_1,\ldots,T_r$ of normal graded matrices over $\ZZ$   is hyperbolic, then $  u_\pm (T_1)+\ldots+ u_\pm (T_r)=0$.

 \section{Cobordism of knots vs. cobordism of matrices}\label{6}

 \subsection{The matrices   of cobordant knots} Setting $R=\ZZ$, we can apply the definitions and results of Section  \ref{5}  to the
graded matrices of knots.   The following  theorem  relates cobordisms of knots to cobordisms of their   matrices.

\begin{theor}\label{bnknk3}  The graded
matrices  of cobordant knots are cobordant.
\end{theor}

The 
proof of this theorem will be given in Section \ref{129ppp8} using the notion of a slice genus introduced in Section \ref{sgenu}. 

\begin{corol}\label{aa439---5}   If two knots $K_1, K_2$ are cobordant, then   	 for any  $m_1,\ldots, m_k \geq 1$, the higher graded matrices $T_\bullet^{m_1,\ldots, m_k }  ( K_1) $ and $T_\bullet^{m_1,\ldots, m_k }  ( K_2)$ are cobordant.
 \end{corol}

\begin{corol}\label{aa439762589317ns5}  The  graded matrices and the higher graded matrices of  slice knots are
hyperbolic.
		   \end{corol}

 \subsection{Slice genus}\label{sgenu}    The {\it   slice  genus}
$sg(K_1,\ldots, K_r)$ of   $r\geq 1$   knots
$K_1\subset \Sigma_1\times \RR$, $\ldots$, $K_r\subset \Sigma_r\times \RR$ is the minimal  integer
$k\geq 0$ satisfying the following condition:
there is an oriented 3-manifold $M$ such that   $ \amalg_{i=1}^r \Sigma_i\subset \partial M$ and  the 1-manifold  $$ \amalg_{i=1}^r K_i\subset   \amalg_{i=1}^r \Sigma_i \times \RR  \subset \partial M \times \RR$$ bounds  in $M \times \RR$ an embedded   compact  (oriented)  surface  of
genus
$  k$.  If there is no such $k$, then   $sg(K_1,\ldots, K_r)=+\infty$. Note that we   do not require
$M$ or the surface in $M\times \RR$ 
 to be connected although it is always possible to achieve this by taking connected sum. The
genus of a disconnected surface is by definition the sum of the genera of its components.

Computing the slice genus of a family of knots is   an interesting geometric problem.  The following lemma estimates the slice genus   from below.

 \begin{lemma}\label{th91373}  For any knots $K_1,\ldots,K_r$,   $$ 
2\, sg(K_1,\ldots,K_r) \geq \sigma
(T_\bullet (K_1),\ldots,T_\bullet (K_r )) .$$
\end{lemma}
\begin{proof}   Consider
a  3-manifold
$M$ as in the definition of the slice genus and  an embedded compact    surface $F\subset M\times \RR$   of genus
$  k=sg(K_1,\ldots,K_r)$  with $\partial F=  \amalg_{i=1}^r K_i\subset     \partial M \times \RR$. 
By assumption,   the components of $\partial F$  project to disjoint
 subsets of $\partial M$. As in the proof of Lemma \ref{l:nblddddd}, we can additionally assume that     $ F\cap (\partial M \times \RR)=\partial F$  and  the restriction to $F$ of the projection $M\times \RR \to M$ is a generic map.  Denote this map $F\to M$ by $\omega$.

The  group   $H_1(F)=H_1(F;\ZZ)$  is generated by the homology classes of $r$ boundary components of $F$
 and a subgroup   $H\subset H_1(F)$ isomorphic to $\ZZ^{2k}$. Set
$$L=
\inc^{-1} (\omega_\ast(H))\subset  H_1(\partial M),$$ where   $\omega_\ast:H_1(F)\to H_1(M)$ is the homomorphism  induced by $\omega$ and $\inc:H_1(\partial M) \to H_1(M)$ is the
inclusion homomorphism.   Denote the homological intersection form $ H_1(\partial M) \times
H_1(\partial M) \to
\ZZ$ by $B$. Since $B$ annihilates the  kernel of
$\inc$,  
\begin{equation}\label{nokl}\rk\,  (B\vert_L:L\times L \to \ZZ) \leq 2 \rk\, \omega_\ast(H) \leq 2\rk  H=4k.\end{equation}

For $t=1,\ldots,r$, consider   the graded matrix $T_t=(G_t,s_t,b_t)$ of the diagram of $K_t$   associated with the projection of $K_i$ to $\partial M$. We derive from   $\omega:F\to M$ a filling $\lambda$ of the tuple $T_1,\ldots , T_r$.  Set $G=\amalg_t G_t$.    As in the proof
of Lemma
\ref{l:nblddddd}, the map
$\omega$ gives rise to an  involution
$\nu$ on the set $$\Join\! \! (\omega\vert_{\partial F}) = \amalg_{t=1}^r (G_t -\{s_t\})= G-\{s_1,\ldots,s_r\}.$$   Each orbit  $X$ of this involution   gives rise to a vector $\lambda_X $ in the lattice $\ZZ G $ freely generated by  
$G$.  This vector is defined as follows. Set $$[X]=\sum_{x\in X} [\omega_x]  \in H_1(\partial M),$$ where   $\omega_x$ is the loop on $\partial M$ associated with the self-crossing $x$ of $\omega\vert_{\partial F}  $ as in the   proof of Lemma
\ref{l:nblddddd}.   This   proof    shows that   $  \inc ([X]) \in   \omega_\ast(H_1(F))  $ or, equivalently,     $$  [X] \in   \inc^{-1}
(\omega_\ast(H_1(F))) .$$  Therefore, there  is a    linear
combination
$$n_X=\sum_{t=1}^r n_{X,t} \, [\omega(K_t)] \in H_1(\partial M)$$ of the homology classes $[\omega(K_1)], \ldots, [\omega(K_r)]
\in H_1(\partial M)$  with $n_{X,t}\in \ZZ$ such that $[X]+n_X\in L$.  (Such $n_X$ may be non-unique; we take any.)  
Set $$\lambda_X= \sum_{x\in X} x +\sum_{t=1}^r n_{X,t}  s_t  \in \ZZ G .$$    The  vectors $\{\lambda_X\}_X$ corresponding to all  orbits $X$  of $\nu$  
 together with  the vector $s_1+\cdots+s_r\in \ZZ G $ form a 
filling   of   the tuple 
$T_1,\ldots,T_r $. The matrix    of this   filling     is obtained by evaluating the form $B$ on the  homology classes
 $$   [X] + n_X  \in H_1(\partial M)\quad {\text {and}}  \quad [\omega(K_1)]+
\ldots+ [\omega(K_r)]\in H_1(\partial M).$$  Since   these homology classes belong to
$L$, Formula \eqref{nokl}   implies that 
  the rank of this matrix 
is smaller than or equal to $4k$. Thus 
$$\sigma (T_1,\ldots, T_r )\leq 2k=
2\, sg(K_1,\ldots,K_r).$$  By definition, $T_\bullet (K_t)=(T_t)_\bullet $ for $t=1,\ldots, r$. 
By Corollary \ref{aaal:gfglpmns5}, $$\sigma (T_1,\ldots, T_r )=\sigma ((T_1)_\bullet, \ldots, (T_r)_\bullet)= \sigma
(T_\bullet (K_1),\ldots,T_\bullet (K_r )).$$
Hence, $ \sigma
(T_\bullet (K_1),\ldots,T_\bullet (K_r ))\leq 2\, sg(K_1,\ldots,K_r)$. 
\end{proof}

\subsection{Proof of Theorem \ref{bnknk3}}\label{129ppp8} If two knots $K_1$, $K_2$ are cobordant, then we have $sg(K_1, -\overline {K}_2)=0$. By Lemma
\ref{th91373}, $\sigma (T_\bullet(K_1), T_\bullet(-\overline {K}_2))=0$.  As we know,  $T_\bullet(-\overline {K}_2)=-T_\bullet(K_2)$. Thus,
 $\sigma (T_\bullet(K_1),  -T_\bullet(K_2))=0$ so that $T_\bullet(K_1)$ is cobordant to ~$T_\bullet(K_2)$.

   \section{Miscellaneous remarks and open questions}\label{7}
   
   Theorem \ref{bnknk3} leads to a number of interesting questions which we   briefly discuss.  
      
    \subsection{Invariants of graded matrices}\label{16--8}   Theorem \ref{bnknk3} suggests to search for a cobordism classification of skew-symmetric graded matrices over $\ZZ$. We list here
    several constructions of cobordism invariants of graded matrices (not necessarily skew-symmetric) which may   help.
    
    As we know, the genus $\sigma(T)$ of a graded matrix $T $ over a domain $R$ is a cobordism invariant of $T$. Moreover, for any finite family of graded matrices $T_1,\ldots, T_r$ of  graded matrices over $R$, the genus $\sigma (T,T_1, \ldots, T_r)$ is a cobordism invariant of~$T$. Further invariants may be obtained by   ring replacements. Given a ring homomorphism $\varphi$ from   $R$ to a domain $R'$, we can derive from $T=(G,s, b)$ the graded matrix $T_\varphi=(G,s, \varphi \circ b)$ over $R'$. The cobordism invariants of $T_\varphi$ will be cobordism invariants of $T$. For example, let $R=\ZZ$, $R'=\ZZ/p\ZZ$, where $p\geq 2$ is a prime integer, and let $\varphi:R\to R'$ be the projection. Then the {\it $p$-genus} $\sigma_p(T)=\sigma (T_\varphi)$ is a   cobordism invariant of  a graded matrix $T$ over $\ZZ$. It is clear that $0\leq \sigma_p (T) \leq \sigma (T)$ and that $\sigma_p(T)=\sigma (T)$ for any given $T$ and all sufficiently big prime integers $p$. Thus, $\sigma(T)=\max_p \sigma_p(T)$. 
    
   More   invariants of graded matrices over a domain $R$ can be obtained using the following transformations of graded matrices. Pick a set $A\subset R$ such that $-A= A$, where $-A=\{-a \,\vert \,  a\in A\}$. For a   graded matrix $T=(G,s,b)$ over $R$, set $G_A=\{g\in G\,\vert \, b(g,s)\in A \}$. The bipartition of $G $ induces a bipartition of $G_A$ by $(G_A)_\pm =G_A\cap G_\pm$. We say that  $T $   is {\it $A$-normal}, if $b(s,s)\in A$ or equivalently, if $s\in G_A$. If $T$ is $A$-normal, then  $\gamma_A(T) =(G_A, s, b\vert_{G_A})$ is a graded matrix over~$R$. 
   It is easy to see that if two $A$-normal graded matrices $T_1, T_2$ are cobordant, then  the graded matrices $\gamma_A(T_1), \gamma_A(T_2)$ are cobordant. Thus,  the formula $T\mapsto \gamma_A(T)$ induces a   transformation on the set of cobordism classes  of $A$-normal graded matrices.  This construction can be applied to a   skew-symmetric graded matrix $T$ over $\ZZ$ and   $A=(-B)\cup \{0\} \cup B$, where $B$ is an arbitrary set of positive integers. 
 In particular,    the genus and the $p$-genera of $\gamma_A(T)$
are cobordism invariants of   $T$ parametrized by $B$.

 \subsection{Realization problem}\label{1---}  It would be interesting to give an algebraic description of all primitive graded matrices over $\ZZ$ that can be realized as the graded matrices of knots. Formula \ref{rede}  yields a  necessary condition.  
 
 \subsection{Secondary obstructions to   sliceness}\label{1hhhgp8} As we know, the graded matrix and the higher graded matrices of knots yield obstructions to  the sliceness of knots. We outline a   Casson-Gordon style construction of secondary obstructions.

Let $K$  be a knot  on a closed connected (oriented) surface   $  \Sigma $ and let $ \underline K$ be the projection of $K$ to $\Sigma$.  Fix an integer $m\geq 2$
and set  
$H= H_1(\Sigma; \ZZ/m\ZZ)$.      For
any $h\in H$ such that $h\cdot \underline K=0$, consider the covering $\Sigma_h \to \Sigma$ corresponding to the normal subgroup of $\pi_1(\Sigma)$ consisting of the homotopy classes of loops on $\Sigma$ whose   intersection number with $h$ is $0$. The knot $K\subset \Sigma \times \RR$ lifts to a knot $K_h \subset \Sigma_h \times \RR$ and the diffeomorphism type of  $K_h$ does not depend on the choice of the lift.    The invariants of   $K_h$ can be viewed as  invariants of $K$ parametrized by $m$ and $h$. In
particular, we can consider the   polynomials $ u_\pm(K_h)$, the graded matrix  $ T_\bullet(K_h)$, the slice genus of $K_h$, etc. In the next lemma, a {\it Lagrangian}  
is an additive  group $L\subset H$   equal to its annihilator  $\Ann (L)=\{g\in H\,\vert \,  g\cdot L =0\}$. %%%If $m$ is   prime, then    each%%%
%%%%MLagrangian $L\subset H$ is a direct summand of $H$ and $H/L\approx L$.%%%

 \begin{lemma}\label{th955555553}  If $K$ is slice, then 
there is a Lagrangian  $L\subset H$ satisfying the following conditions: $[\underline K]\in L$;  the knot $K_h$ is
slice for all $h\in L$;   for any diagram $D$ of $K$ on $\Sigma$,   there is an involution
  on the set $ \Join \!\! (D) $ such that  for any  its orbit $X\subset\, \, \Join \!\! (D) $,  we have  $\sum_{x\in X} [D_x] \in L$.
\end{lemma}

 This lemma implies that for  all $h\in L$,    the graded matrix $T_\bullet (K_h)$ is hyperbolic.
 
\begin{proof}   If $K$ is slice, then there  is an oriented 3-manifold $M$ and an embedding  $ \Sigma \subset \partial M$ such that  the image of $K $ under the induced embedding
 $  \Sigma \times \RR  \subset \partial M \times \RR$ bounds an embedded disk   $B\subset M \times \RR$. Replacing if necessary $M$ by a compact submanifold containg both $\Sigma$ and the projection of $B$ to $M$, we can assume that $M$ is compact. Gluing handlebodies to all the components of $\partial M$ distinct from $\Sigma$, we can additionally assume that $\partial M=\Sigma$.  
 
Set $R=\ZZ/m\ZZ$. Consider  the boundary  homomorphism $\partial : H_2(M,\partial M;R) \to H_1(\partial M;R)  =H$ and  the inclusion
homomorphism
$i:H=H_1(\partial M;R) \to H_1(M; R)$. Set $L=\Ima (\partial)=\Ker (i)$. 
  It is well known that
  $L  $ is a Lagrangian.  For completeness, we outline a proof. An element $g\in H$ belongs to $ \Ann (L)$ iff  $g \cdot  \partial a =0$ for every 
$a\in H_2(M,\partial M;R)$. Given $a$,   the Poincar\'e duality says that there is a unique $\widetilde a \in H^1(M;R)=\Hom (H_1(M;R), R)$ such that $a=\widetilde a \cap
[M]$.  Then  
$g\cdot  \partial a = i(g) \cdot a =
\widetilde a (i(g))$.   Therefore $g\in \Ann (L)$ iff  $i(g)$ is annihilated by all homomorphisms $H_1(M;R) \to R$. 
This holds iff
$i(g)=0$, that is iff
$g\in L$. 

Pick  
$h\in L$ and  pick   any $a $ in  $\partial^{-1} (h)\subset  H_2(M,\partial M;R)$. The   cohomology class $\widetilde
a\in H^1(M;R)$ defines a cyclic
$m$-fold covering  
$\widetilde M \to M$. The disk $B\subset M \times \RR$ lifts   to a disk $\widetilde  B\subset \widetilde M \times \RR$. By  
$\widetilde a (i(g))=  g\cdot \partial a= g\cdot h$ for all $ g\in H$,  the covering  $\widetilde M \to M$ restricted to $\Sigma=\partial M$ can be identified with    the covering $\Sigma_h \to \Sigma$ considered above. Then $K_h=\partial \widetilde B \subset \Sigma_h  \times \RR$ is a    slice knot. Constructing an involution  
 on  $ \Join \!\! (D) $ as in the proof of Lemma \ref{l:nblddddd}  (for $F=B$) we obtain the last claim of the lemma.
\end{proof}

 \subsection{More on the slice genus}\label{16668}  Consider again the slice  genus 
$sg(K_1,\ldots, K_r)\geq 0$ of a family  of     knots
$K_1 $, $\ldots$, $K_r $.     The slice genus  does not depend on the order in the sequence $K_1,\ldots, K_r$ and  is preserved if  
$K_1,\ldots, K_r$ are replaced  with  cobordant knots.       If  $sg(K_1,\ldots, K_r)=0$ then we call the family of knots 
$K_1,\ldots, K_r$  {\it slice}.  One immediate corollary of 
Lemma \ref{th91373}
is the following theorem.

\begin{theor}\label{tqpaw53}  If a finite family  of knots    is slice, then   the family of their   graded
matrices is  hyperbolic.
\end{theor}

 Observe, in generalization  of Lemma \ref{le},  that 
 if the family $K_1,\ldots, K_r$ is slice, then  the family of the  $m$-th coverings 
$   K_1^{(m)},\ldots,   K_r^{(m)}$ is slice for any   $m\geq 1$.
By Theorem \ref{tqpaw53}, if $K_1,\ldots, K_r$ is slice, then for any finite sequence of   integers $m_1,\ldots,m_k\geq 1$, the family of graded matrices  $
T_\bullet^{(m_1,\ldots,m_k)}( K_1),\ldots ,  T_\bullet^{(m_1,\ldots,m_k)}( K_r) $ is hyperbolic.

For $r=2$, we can rewrite the slice genus    in the equivalent form
$sg'(K_1, K_2)= sg(K_1, -K_2^-)$.    The
  number  $sg'(K_1, K_2)$ depends only on the cobordism classes of
$K_1, K_2$ and     defines a metric on the set of cobordism classes of knots.  

In \cite{tu2}  the author defined a slice genus $sg $
for any finite family of disjoint loops   on surfaces.   It follows from the definitions that   $sg(K_1,\ldots, K_r) \geq sg(\underline {K}_1,\ldots, \underline {K}_r)$.

  \subsection{Long knots}\label{rems2} The definitions and results   of this paper can be generalized  to links and to long knots on surfaces. The context of long knots is especially interesting since their cobordism  classes form a group. It is a challenging question to compute this group.   Is it abelian ?

                     \end{document}